# Fréchet differentiability of the metric projection operator in Banach spaces

Jinlu Li

Department of Mathematics
Shawnee State University
Portsmouth, Ohio 45662
USA

**Abstract**

In this paper, we prove Fréchet differentiability of the metric projection operator onto closed balls, closed and convex cylinders and positives cones in uniformly convex and uniformly smooth Banach spaces. With respect to these closed and convex subsets, we find the exact expressions for Fréchet derivatives and Gâteaux directional derivatives of the metric projection operator.

1. **Introduction and preliminaries**

Throughout this paper, unless otherwise stated, let $(X, \|\cdot\|)$ be a real uniformly convex and uniformly smooth Banach space with topological dual space $(X^*, \|\cdot\|_*)$. Let $\langle \cdot, \cdot \rangle$ denote the real canonical pairing between $X^*$ and $X$. Let $C$ be a nonempty closed and convex subset of $X$. Let $P_C: X \to C$ denote the (standard) metric projection operator, which is a well-defined single-valued mapping. In operator theory, it is clear that the metric projection operator $P_C$ is one of the most important operators, which has been studied by many authors with a long history (see [1, 2, 5, 8, 12, 15, 17, 21]). The metric projection operator $P_C$ has been widely applied to nonlinear analysis, optimization theory, approximation theory, fixed point theory, variational inequalities, so and forth (see [1, 5, 15]).

In operator theory, one of the most important research topics is the continuity of the considered operators. In addition to continuity, the smoothness of operators has attracted a lot of attention in the mathematics community. Several types of differentiability of operators, in particular, the metric projection operator, have been introduced and studied in Hilbert spaces (see [6, 7, 9, 10, 13, 14, 18]); in Banach spaces and normed linear spaces (see [4, 7, 11, 16, 19, 20]).

In this paper, we focus on the differentiability of the metric projection operator in uniformly convex and uniformly smooth Banach spaces. We recall some concepts of differentiability of the metric projection operator in uniformly convex and uniformly smooth Banach space $X$ below, which will be used in this paper.

I. Gâteaux directional differentiability of the metric projection operator in uniformly convex and uniformly smooth Banach spaces (see Definition 4.1 in [11]). For $x \in X$ and $v \in X$ with $v \neq \theta$, if the following limit exists,

$$P'_C(x)(v) := \lim_{t \downarrow 0} \frac{P_C(x+tv) - P_C(x)}{t},$$

then, $P_C$ is said to be Gâteaux directionally differentiable at point $x$ along direction $v$. $P'_C(x)(v)$ is called the Gâteaux directional derivative of $P_C$ at point $x$ along direction $v$; and $v$ is called a Gâteaux differentiable direction of $P_C$ at $x$.

II. Fréchet differentiability of the metric projection operator in uniformly convex and uniformly smooth Banach spaces (see Definition 1.13 in [17]). For any given $\bar{x} \in X$, if there is a linear and continuous mapping $\nabla P_C(\bar{x}): X \to X$ such that

$$\lim_{x \to \bar{x}} \frac{P_C(x) - P_C(\bar{x}) - \nabla P_C(\bar{x})(x-\bar{x})}{\|x - \bar{x}\|} = \theta.$$

then $P_C$ is said to be Fréchet differentiable at $\bar{x}$ and $\nabla P_C(\bar{x})$ is called the Fréchet derivative of $P_C$ at $\bar{x}$.

III. Strict Fréchet differentiability of the metric projection operator in uniformly convex and uniformly smooth Banach spaces (see Definition 1.13 in [17]). For any given $\bar{x} \in X$, if there is a linear and continuous mapping $\nabla P_C(\bar{x}): X \to X$ such that

$$\lim_{(u,v) \to (\bar{x},\bar{x})} \frac{P_C(u) - P_C(v) - \nabla P_C(\bar{x})(u-v)}{\|u-v\|} = \lim_{u \to \bar{x}, v \to \bar{x}} \frac{P_C(u) - P_C(v) - \nabla P_C(\bar{x})(u-v)}{\|u-v\|} = \theta,$$

then $P_C$ is said to be strictly Fréchet differentiable at $\bar{x}$ at $\bar{x}$ and $\nabla P_C(\bar{x})$ is called the Fréchet derivative of $P_C$ at $\bar{x}$.

In addition to the above three concepts of differentiability of the metric projection operator, in [17], Mordukhovich introduced the concept of generalized differentiation of operators and provide many useful properties, such as differential calculus in Banach spaces. The theory of generalized differentiation of operators in Banach spaces has been widely applied to nonlinear analysis, which are presented in this book [17].

The above three types of differentiability have the following inclusion properties (see [13] for mor details): for any given point $\bar{x} \in X$,

$P_C$ is strictly Fréchet differentiable at $\bar{x}$

$\Rightarrow P_C$ is Fréchet differentiable at $\bar{x}$

$\Rightarrow P_C$ is Gâteaux directionally differentiable at $\bar{x}$.

However, the converse statements may not hold. For some point $\bar{x} \in X$, we have

$P_C$ is Gâteaux directionally differentiable at $\bar{x}$

$\not\Rightarrow$ $P_C$ is Fréchet differentiable at this point $\bar{x}$.

For the study of the differentiability of the metric projection operator in Banach spaces, we started at the simpler case that is the Gâteaux directional differentiability. In [11], the Gâteaux directional differentiability of the metric projection operator in uniformly convex and uniformly smooth Banach spaces was studied, in which, some properties of Gâteaux directional derivatives of $P_C$ are proved. Moreover, when the considered subset $C$ is a closed ball or the positive cone of the considered Banach space, the exact expressions for Gâteaux directional derivatives of $P_C$ are presented in this paper [11]. Following the ideas in [11], the Gâteaux directional differentiability of the metric projection operator was studied in Hilbert spaces in [14].

Very recently, in [13], the present author studied the strict Fréchet differentiability of the metric projection operator in Hilbert spaces. Going a step further, when the considered subset $C$ is a ball, or the positive cone in Euclidean spaces and in the real Hilbert space $l_2$, the precise solutions for the Fréchet derivatives of $P_C$ are provided in [13].

In this paper, we study the Fréchet differentiability of the metric projection operator in uniformly convex and uniformly smooth Banach spaces. We investigate the properties and solutions of the Fréchet derivatives of $P_C$.

In section 2, we prove some properties of the normalized duality mapping and the metric projection operator in uniformly convex and uniformly smooth Banach spaces; in section 3, we prove the Fréchet differentiability of the metric projection operator onto closed balls in uniformly convex and uniformly smooth Banach spaces and we find the exact solutions of the Fréchet derivatives; in section 4, we prove the Fréchet differentiability of the metric projection operator onto closed and convex cylinders in the real Banach space $l_p$, for some $p$ satisfying $1 < p < \infty$ and we find the exact solutions of the Fréchet derivatives; In section 5, we disprove the Fréchet differentiability of the metric projection operator onto the positive cone in real Banach space $L_p(S)$ for some $p$ satisfying $1 < p < \infty$.

## 2. The normalized duality mapping and the metric projection operator in uniformly convex and uniformly smooth Banach spaces

Let $(X, \|\cdot\|)$ be a real uniformly convex and uniformly smooth Banach space with topological dual space $(X^*, \|\cdot\|_*)$. Let $\langle \cdot, \cdot \rangle$ denote the real canonical pairing between $X^*$ and $X$. Let $\theta$ and $\theta^*$ denote the origins and $\mathbb{B}$ and $\mathbb{B}^*$ denote the unit closed balls in $X$ and $X^*$, respectively. It follows that, for any $r > 0$, $r\mathbb{B}$ and $r\mathbb{B}^*$ are closed balls with radius $r$ and centered at the origins in $X$ and $X^*$, respectively. Let $\mathbb{S}$ be the unit sphere in $X$. Then, $r\mathbb{S}$ is the sphere in $X$ with radius $r$ and center $\theta$. For any $c \in X$ and $r > 0$, let $\mathbb{B}(c, r)$ denote the closed ball in $X$ with radius $r$ and center $c$. It follows that $\mathbb{B}(\theta, 1) = \mathbb{B}$ and $\mathbb{B}(\theta, r) = r\mathbb{B}$. Let $I_X$ denote the identity mapping in $X$. The dual mapping of $I_X$ is denoted by $I_X^*$, which is the identity mapping in $X^*$. Let $\mathbb{R}$ denote the set of real numbers.

### 2.1. The normalized duality mapping

Let $J: X \to X^*$ and $J^*: X^* \to X$ be the normalized duality mappings. We have

(i) $\langle J(x), x \rangle = \|x\| \|J(x)\|_* = \|x\|^2 = \|J(x)\|_*^2$, for any $x \in X$;

(ii) $\langle x^*, J^*(x^*) \rangle = \|J^*(x^*)\| \|x^*\|_* = \|J^*(x^*)\|^2 = \|x^*\|_*^2$, for any $x^* \in X^*$;

(iii) In uniformly convex and uniformly smooth Banach spaces, in general,

$$\langle J(x), y \rangle = 0 \quad \not\Rightarrow \quad \langle J(y), x \rangle = 0, \text{ for } x, y \in X.$$

The norm of uniformly convex and uniformly smooth Banach space $X$ is Gâteaux directionally differentiable at every point (see [1, 2, 15, 21]). That is, the following limit exists

$$\lim_{t\downarrow 0}\frac{\|x+ty\|-\|x\|}{t}, \quad \text{uniformly for } (x,y) \in \mathbb{S}\times\mathbb{S}. \tag{2.1}$$

The normalized duality mapping $J$ of the considered uniformly convex and uniformly smooth Banach space $X$ has the following useful properties.

**Lemma 2.1**. *Let $X$ be a uniformly convex and uniformly smooth Banach space. For any $x, y \in X$, one has*

$$2\langle J(y), x-y\rangle \leq \|x\|^2 - \|y\|^2 \leq 2\langle J(x), x-y\rangle.$$

*This implies*

$$\left|\|x\|^2 - \|y\|^2\right| \leq 2(|\langle J(x), x-y\rangle| + |\langle J(y), x-y\rangle|)$$

*Proof.* The proof of the first inequality can be found in Takahashi [21]. We similarly prove the second part.

$$\|x\|^2 - \|y\|^2 - 2\langle J(x), x-y\rangle$$
$$= -\|x\|^2 - \|y\|^2 + 2\langle J(x), y\rangle$$
$$\leq -\|x\|^2 - \|y\|^2 + 2\|x\|\|y\|$$
$$= -(\|x\|-\|y\|)^2 \leq 0. \qquad \square$$

For the considered uniformly convex and uniformly smooth Banach space $X$, we define a function $\Psi: X\times X \to \mathbb{R}$ by

$$\Psi(x,y) = \lim_{t\downarrow 0}\frac{\|x+ty\|-\|x\|}{t}, \quad \text{for any } (x,y) \in X\times X.$$

$\Psi$ is called the function of smoothness of $X$. By property (xi), we can check that $\Psi$ is well-defined satisfying

$$\Psi(\theta, y) = \|y\|, \text{ for any } y \in X.$$

When $(x,y) \in X\times X$, $\Psi(x,y)$ is written as $\psi(x,y)$.

**Lemma 2.2**. *For any $x, y \in X$ with $x \neq \theta$, one has*

$$\Psi(x,y) = \lim_{t\downarrow 0}\frac{\|x+ty\|-\|x\|}{t} = \frac{\langle J(x), y\rangle}{\|x\|}.$$

*In particular, for $(x,y) \in \mathbb{S}\times\mathbb{S}$, by (2.1), one has*

$$\psi(x,y) = \lim_{t\downarrow 0}\frac{\|x+ty\|-\|x\|}{t} = \langle J(x), y\rangle, \quad \text{uniformly for } (x,y) \in \mathbb{S}\times\mathbb{S}. \tag{2.2}$$

*Proof.* For any $x, y \in X$ with $x \neq \theta$, by the first inequality in Lemma 2.1, we have

$$\lim_{t \downarrow 0} \frac{\|x+ty\| - \|x\|}{t} = \lim_{t \downarrow 0} \frac{\|x+ty\|^2 - \|x\|^2}{t(\|x+ty\|+\|x\|)}$$

$$\geq \lim_{t \downarrow 0} \frac{2\langle J(x), (x+ty)-x \rangle}{t(\|x+ty\|+\|x\|)}$$

$$= \lim_{t \downarrow 0} \frac{2\langle J(x), ty \rangle}{t(\|x+ty\|+\|x\|)}$$

$$= \lim_{t \downarrow 0} \frac{2\langle J(x), y \rangle}{\|x+ty\|+\|x\|}$$

$$= \frac{\langle J(x), y \rangle}{\|x\|}.$$

On the other hand, by the second inequality in Lemma 2.1 and by the continuity of the normalized duality mapping $J$, we have

$$\lim_{t \downarrow 0} \frac{\|x+ty\| - \|x\|}{t} = \lim_{t \downarrow 0} \frac{\|x+ty\|^2 - \|x\|^2}{t(\|x+ty\|+\|x\|)}$$

$$\leq \lim_{t \downarrow 0} \frac{2\langle J(x+ty), (x+ty)-x \rangle}{t(\|x+ty\|+\|x\|)}$$

$$= \lim_{t \downarrow 0} \frac{2\langle J(x+ty), ty \rangle}{t(\|x+ty\|+\|x\|)}$$

$$= \lim_{t \downarrow 0} \frac{2\langle J(x+ty), y \rangle}{\|x+ty\|+\|x\|}$$

$$= \frac{\langle J(x), y \rangle}{\|x\|}. \qquad \square$$

### 2.2. The metric projection operator

Let $C$ be a nonempty closed and convex subset of this uniformly convex and uniformly smooth Banach space $X$. Let $P_C \colon X \to C$ denote the (standard) metric projection operator. For any $x \in X$, $P_C x \in C$ such that

$$\|x - P_C x\| \leq \|x - z\|, \text{ for all } z \in C.$$

The annihilator of $C$ is denoted by $C^\perp$, which is a subspace of $X^*$ and is defined by

$$C^\perp = \{ g \in X^* \colon \langle g, x \rangle = 0, \text{ for all } x \in C \}.$$

The basic variational principle of the metric projection operator $P_C$ is, for any $x \in X$ and $u \in C$,

$$u = P_C x \quad \Leftrightarrow \quad \langle J(x - u), u - z \rangle \geq 0, \text{ for all } z \in C,$$

In particular, if $M$ is a closed subspace of $X$, then, the basic variational principle of the metric projection operator $P_M \colon X \to M$ becomes to the following equality version.

**Lemma 2.3.** *Let M be a closed subspace of X. For any $x \in X$ and $u \in M$,*

$$u = P_M x \quad \Longleftrightarrow \quad \langle J(x-u), u-z \rangle = 0, \text{ for all } z \in M,$$

$$\Longleftrightarrow \quad \langle J(x-u), z \rangle = 0, \text{ for all } z \in M,$$

$$\Longleftrightarrow \quad J(x-u) \in M^\perp.$$

The metric projection operator $P_M: X \to M$ is conditionally linear, which is demonstrated by the following lemma (see Lemma 1.4 in Alber [1, 2]).

**Lemma 2.4.** *Let M be a closed subspace of X. For arbitrary $x \in X$ and $y \in M$, we have*

$$P_M(\alpha x + \beta y) = \alpha P_M(x) + \beta y, \text{ for any } \alpha, \beta \in \mathbb{R}.$$

*In particular, by taking $y = \theta$, we have*

$$P_M(\alpha x) = \alpha P_M(x), \text{ for any } \alpha \in \mathbb{R}.$$

*Proof.* Notice that $\alpha P_M(x) + \beta y \in M$, for any $\alpha, \beta \in \mathbb{R}$. By the properties of the normalized duality mapping $J$, we have

$$\langle J((\alpha x + \beta y) - (\alpha P_M(x) + \beta y)), z \rangle$$

$$= \langle J(\alpha x - \alpha P_M(x)), z \rangle$$

$$= \langle J(\alpha(x - P_M(x))), z \rangle$$

$$= \alpha \langle J(x - P_M(x)), z \rangle$$

$$= 0, \text{ for all } z \in M.$$

By Lemma 2.3, this proves Lemma 2.4. □

For any given $\bar{x} \in X \setminus \{\theta\}$, let $S(\bar{x})$ denote the one-dimensional (closed) subspace of $X$ generated by this vector $\bar{x}$. Then, the annihilator of $S(\bar{x})$ satisfies

$$S(\bar{x})^\perp = \{g \in X^*: \langle g, \bar{x} \rangle = 0\}.$$

**Proposition 2.5.** *For any given $\bar{x} \in X \setminus \{\theta\}$, the metric projection operator $P_{S(\bar{x})}: X \to S(\bar{x})$ has the following properties. For any $x \in X$, one has*

(i) $P_{S(\bar{x})}(\alpha x + \beta \bar{x}) = \alpha P_{S(\bar{x})}(x) + \beta \bar{x}, \text{ for any } \alpha, \beta \in \mathbb{R}.$

*In particular,*

$$P_{S(\bar{x})}(\alpha x) = \alpha P_{S(\bar{x})}(x), \text{ for any } \alpha \in \mathbb{R};$$

(ii) $\langle J(x - P_{S(\bar{x})}(x)), \bar{x} \rangle = 0, \text{ for all } x \in X;$

(iii) *If $X$ is a Hilbert space, then*

$$P_{S(\bar{x})}(x) = \frac{\langle \bar{x}, x \rangle}{\|\bar{x}\|^2} \bar{x}, \text{ for any } x \in X;$$

(iv) *If $X$ is not a Hilbert space, then, in general*

$$P_{S(\bar{x})}(x) \neq \frac{\langle J(\bar{x}), x \rangle}{\|\bar{x}\|^2} \bar{x}, \text{ for } x \in X.$$

*Proof.* The proofs of parts (i–iii) are straight forward and they are omitted here. We only prove (iv). We prove (iv) by constructing a counter example. By Lemma 2.4 and part (ii) of this lemma, for any $x \in X$ and $u \in S(\bar{x})$, we have

$$u = P_{S(\bar{x})}(x) \quad \Leftrightarrow \quad \langle J(x-u), \bar{x} \rangle = 0. \tag{2.3}$$

Now, we provide a counter example for proving part (iv). Let $X = \mathbb{R}^3$ equipped with the 3-norm $\|\cdot\|_3$ defined by, for any $z = (z_1, z_2, z_3) \in \mathbb{R}^3$,

$$\|z\|_3 = \sqrt[3]{|z_1|^3 + |z_2|^3 + |z_3|^3}.$$

Then, $(\mathbb{R}^3, \|\cdot\|_3)$ is a uniformly convex and uniformly smooth Banach space that is not a Hilbert space. The dual space of $(\mathbb{R}^3, \|\cdot\|_3)$ is $(\mathbb{R}^3, \|\cdot\|_{\frac{3}{2}})$. The normalized duality mapping $J: \mathbb{R}^3 \to \mathbb{R}^3$ satisfies the following conditions, for any $z = (z_1, z_2, z_3) \in \mathbb{R}^3$ with $z \neq 0$,

$$Jz = \left( \frac{|z_1|^2 \text{sign}(z_1)}{\|z\|_3}, \frac{|z_2|^2 \text{sign}(z_2)}{\|z\|_3}, \frac{|z_3|^2 \text{sign}(z_3)}{\|z\|_3} \right).$$

Take $\bar{x} = (1, 2, 1)$ and $x = (3, -2, 1)$. Then, we have $\|\bar{x}\|_3 = \sqrt[3]{10}$ and

$$J(\bar{x}) = \left( \frac{1}{\sqrt[3]{10}}, \frac{4}{\sqrt[3]{10}}, \frac{1}{\sqrt[3]{10}} \right).$$

This implies

$$\frac{\langle J(\bar{x}), x \rangle}{\|\bar{x}\|^2} \bar{x} = \frac{\frac{-4}{\sqrt[3]{10}}}{(\sqrt[3]{10})^2} \bar{x} = \left( -\frac{2}{5}, -\frac{4}{5}, -\frac{2}{5} \right).$$

It follows that

$$x - \frac{\langle J(\bar{x}), x \rangle}{\|\bar{x}\|^2} \bar{x} = (3, -2, 1) - \left( -\frac{2}{5}, -\frac{4}{5}, -\frac{2}{5} \right) = \left( \frac{17}{5}, -\frac{6}{5}, \frac{7}{5} \right).$$

And

$$J\left( x - \frac{\langle J(\bar{x}), x \rangle}{\|\bar{x}\|^2} \bar{x} \right) = J\left( \frac{17}{5}, -\frac{6}{5}, \frac{7}{5} \right) = \frac{\frac{(17^2, -6^2, 7^2)}{5^2}}{\sqrt[3]{\frac{17^3 + 6^3 + 7^3}{5}}} = \frac{(17^2, -6^2, 7^2)}{5\sqrt[3]{17^3 + 6^3 + 7^3}}.$$

We calculate

$$\langle J\left( x - \frac{\langle J(\bar{x}), x \rangle}{\|\bar{x}\|^2} \bar{x} \right), \bar{x} \rangle = \langle \frac{(17^2, -6^2, 7^2)}{5\sqrt[3]{17^3 + 6^3 + 7^3}}, (1, 2, 1) \rangle > 0.$$

By (2.3), this implies

$$P_{S(\bar{x})}(x) \neq \frac{\langle J(\bar{x}), x \rangle}{\|\bar{x}\|^2} \bar{x}. \qquad \square$$

**Lemma 2.6.** *Let $X$ be a uniformly convex and uniformly smooth Banach space. Let $\bar{x} \in X$ with $\bar{x}$*

$\neq \theta$. *Then the following subset of X*

$$\{x \in X: \langle Jx, \bar{x} \rangle = 0\}$$

*is a closed cone in X with vertex* $\theta$. *However, in general, it is not convex.*

*Proof.* Since the normalized duality mapping $J$ is a continuous and homogeneous operator, it follows that the subset $\{x \in X: \langle Jx, \bar{x} \rangle = 0\}$ is a closed cone with vertex at $\theta$ in $X$. Next, we construct a counter-example to show that, the set $\{x \in X: \langle Jx, \bar{x} \rangle = 0\}$ is not convex, in general.

Let $X = \mathbb{R}^3$ be given in Proposition 2.5. Take $\bar{x} = (25, 37, 77)$, $v = (3, -2, -1)$, $w = (1, -3, 2)$. Then $\|v\|_3 = \|w\|_3$. By (2.3), we have

$$Jv = \left(\frac{9}{\|v\|_3}, \frac{-4}{\|v\|_3}, \frac{-1}{\|v\|_3}\right) \quad \text{and} \quad Jw = \left(\frac{1}{\|w\|_3}, \frac{-9}{\|w\|_3}, \frac{4}{\|w\|_3}\right).$$

Then, we have

$$\langle Jv, \bar{x} \rangle = 0 \quad \text{and} \quad \langle Jw, \bar{x} \rangle = 0.$$

It follows that $v, w \in \{x \in X: \langle Jx, \bar{x} \rangle = 0\}$. Take a convex combination of $v$ and $w$ by

$$g = \frac{2}{3}v + \frac{1}{3}w = \left(\frac{7}{3}, -\frac{7}{3}, 0\right).$$

It follows that $\|g\|_3 = \frac{7}{3}\sqrt[3]{2}$. We calculate

$$Jg = \left(\frac{\left|\frac{7}{3}\right|^2 \text{sign}(\frac{7}{3})}{\frac{7}{3}\sqrt[3]{2}}, \frac{\left|-\frac{7}{3}\right|^2 \text{sign}(-\frac{7}{3})}{\frac{7}{3}\sqrt[3]{2}}, \frac{|0|^2 \text{sign}(0)}{\frac{7}{3}\sqrt[3]{2}}\right)$$

$$= \frac{7}{3\sqrt[3]{2}}(1, -1, 0).$$

This implies

$$\langle Jg, \bar{x} \rangle = \langle \frac{7}{3\sqrt[3]{2}}(1, -1, 0), (25, 37, 77) \rangle = -14\sqrt[3]{4} < 0.$$

This shows that $\frac{2}{3}v + \frac{1}{3}w \notin \{x \in X: \langle Jx, \bar{x} \rangle = 0\}$; and therefore, the set $\{x \in X: \langle Jx, \bar{x} \rangle = 0\}$ is not convex. □

**Corollary 2.7.** *For any given $\bar{x} \in X \setminus \{\theta\}$, the metric projection operator $P_{S(\bar{x})}: X \to S(\bar{x})$ has the following properties*

(a) $P_{S(\bar{x})}: X \to S(\bar{x})$ *is not linear;*
(b) *For $y \in S(\bar{x}) \setminus \{\theta\}$, the inverse image of $y$ with $P_{S(\bar{x})}$*

$$P_{S(\bar{x})}^{-1}(y) \text{ is not convex, in general.}$$

*Proof.* We use the results in Lemma 2.6 to construct a counter example to prove both (a) and (b) of this corollary. For $x, y \in X$, by Lemma 2.4, we have

$$\langle J(x - P_{S(\bar{x})}(x)), \bar{x} \rangle = 0 \quad \text{and} \quad \langle J(y - P_{S(\bar{x})}(y)), \bar{x} \rangle = 0.$$

Let $X = \mathbb{R}^3$ and let $\bar{x} = (25, 37, 77)$, $v = (3, -2, -1)$ and $w = (1, -3, 2)$ be given in Lemma 2.6. Let $x = v + \bar{x}$ and $y = w + \bar{x}$. From the proof of Lemma 2.6, we have

$$\langle J(x - \bar{x}), \bar{x} \rangle = \langle J(v), \bar{x} \rangle = 0 \quad \text{and} \quad \langle J(y - \bar{x}), \bar{x} \rangle = \langle J(w), \bar{x} \rangle = 0$$

By Lemma 2.4, this implies

$$P_{S(\bar{x})}(x) = \bar{x} \quad \text{and} \quad P_{S(\bar{x})}(y) = \bar{x}.$$

Let $g = \frac{2}{3}v + \frac{1}{3}w$ be given in the proof of Lemma 2.6. Take a convex combination of $x$ and $y$ by

$$z = \frac{2}{3}x + \frac{1}{3}y = \frac{2}{3}v + \frac{1}{3}w + \bar{x} = g + \bar{x}.$$

By the proof of Lemma 2.6, we have

$$\langle J(z - \bar{x}), \bar{x} \rangle = \langle J(g), \bar{x} \rangle = -14\sqrt[3]{4} < 0.$$

By Lemma 2.4 again, this implies

$$P_{S(\bar{x})}(z) \neq \bar{x}.$$

It follows that $P_{S(\bar{x})}^{-1}(\bar{x})$ is not convex. □

In particular, when the metric projection is onto a closed ball in $X$, it has the following analytic representations. We list it below as a lemma without proof.

**Lemma 2.8.** *For any $r > 0$, the metric projection $P_{r\mathbb{B}}: X \to r\mathbb{B}$ satisfies the following formula.*

$$P_{r\mathbb{B}}(x) = \begin{cases} x, & \text{for any } x \in r\mathbb{B}, \\ \frac{r}{\|x\|}x, & \text{for any } x \notin r\mathbb{B}. \end{cases} \tag{2.4}$$

### 2.3. Some differentiability properties of the metric projection operator

Let $C$ be a nonempty closed and convex subset of $X$. Similarly, to Hilbert cases studied in [13], we consider some properties of the differentiability of $P_C$ in uniformly convex and uniformly smooth Banach spaces.

**Proposition 2.9.** *Let $C$ be a nonempty closed and convex subset of $X$. Let $y \in C$. Suppose $(P_C^{-1}(y))^o \neq \emptyset$. Then, $P_C$ is strictly Fréchet differentiable on $(P_C^{-1}(y))^o$ such that,*

$$\nabla P_C(\bar{x}) = \theta, \text{ for any } \bar{x} \in (P_C^{-1}(y))^o.$$

*That is, for any $\bar{x} \in (P_C^{-1}(y))^o$, one has*

$$\nabla P_C(\bar{x})(w) = \theta, \text{ for every } w \in X \text{ with } w \neq \theta.$$

*Proof.* The proof is same with the proof of Proposition 2.2 in [13]. It is omitted here. □

**Proposition 2.10.** *Let C be a nonempty closed and convex subset of X. Suppose $C^o \neq \emptyset$. Then $P_C$ is strictly Fréchet differentiable on $C^o$ satisfying*

$$\nabla P_C(\bar{x}) = I_X, \text{ for any } \bar{x} \in C^o.$$

*That is, for any $\bar{x} \in C^o$, one has*

$$\nabla P_C(\bar{x})(w) = w, \text{ for every } w \in X \text{ with } w \neq \theta.$$

*Proof.* The proof is same with the proof of Proposition 2.3 in [13]. It is omitted here. □

### 2.4. Semi-orthogonal decompositions

With respect to an arbitrarily given $\bar{x} \in X \backslash \{\theta\}$, every $x \in X$ can be written as

$$x = \frac{\langle J(\bar{x}), x \rangle}{\|\bar{x}\|^2} \bar{x} + \left( x - \frac{\langle J(\bar{x}), x \rangle}{\|\bar{x}\|^2} \bar{x} \right), \text{ for all } x \in X. \tag{2.5}$$

It is clear that

$$\langle J(\bar{x}), x - \frac{\langle J(\bar{x}), x \rangle}{\|\bar{x}\|^2} \bar{x} \rangle = 0, \text{ for all } x \in X.$$

Hence, (2.5) is called a semi-orthogonal decomposition of $x \in X$ with respect to the arbitrarily given $\bar{x} \in X \backslash \{\theta\}$. Since $\frac{\langle J(\bar{x}), x \rangle}{\|\bar{x}\|^2} \bar{x} \in S(\bar{x})$, for all $x \in X$, we define a real valued functional $a(\bar{x}; \cdot): X \to \mathbb{R}$ by

$$a(\bar{x}; x) = \frac{\langle J(\bar{x}), x \rangle}{\|\bar{x}\|^2}, \text{ for all } x \in X.$$

Let

$$O(\bar{x}) := \{x \in X : \langle J(\bar{x}), x \rangle = 0\}$$

We define a mapping $o(\bar{x}; \cdot): X \to O(\bar{x})$ by

$$o(\bar{x}; x) := x - \frac{\langle J(\bar{x}), x \rangle}{\|\bar{x}\|^2} \bar{x} = x - a(\bar{x}; x)\bar{x}, \text{ for all } x \in X.$$

The following lemma provides some properties of $a(\bar{x}; \cdot)$ and $o(\bar{x}; \cdot)$. These properties will play important roles and will be repeatedly used in the following section of this paper.

**Proposition 2.11.** *For any fixed $\bar{x} \in X \backslash \{\theta\}$, the real valued function $a(\bar{x}; \cdot)$ and the mapping $o(\bar{x}; \cdot)$ have the following properties.*

(i) $a(\bar{x}; \cdot): X \to \mathbb{R}$ *is a real valued linear and continuous function satisfying*

$$\langle J(\bar{x}), x \rangle = a(\bar{x}; x)\|\bar{x}\|^2, \quad \text{for any } x \in X;$$

(ii) $o(\bar{x}; \cdot): X \to O(\bar{x})$ *is a linear and continuous mapping satisfying*

$$\langle J(\bar{x}), o(\bar{x}; x) \rangle = 0, \quad \text{for any } x \in X;$$

(iii) $u \to \bar{x} \iff a(\bar{x}; u) \to 1$ and $o(\bar{x}; u) \to \theta$, for $u \in X$;

(iv) $\|u\| \geq \|a(\bar{x}; u)\bar{x}\|$, for any $u \in X$;

(v) $\|u\| \geq \frac{1}{2}\|o(\bar{x}; u)\|$, for any $u \in X$.

*Proof.* By (2.5), the proofs of (i, ii, iii) are straight forward and they are omitted here. We prove (iv). For any $u \in X$, by part (ii), we have

$$\|a(\bar{x}; u)\bar{x}\|\|u\|$$
$$= \|J(a(\bar{x}; u)\bar{x})\|_*\|u\|$$
$$\geq |\langle J(a(\bar{x}; u)\bar{x}), u\rangle|$$
$$= |\langle J(a(\bar{x}; u)\bar{x}), a(\bar{x}; u)\bar{x} + o(\bar{x}; u)\rangle|$$
$$= |\langle J(a(\bar{x}; u)\bar{x}), a(\bar{x}; u)\bar{x}\rangle + a(\bar{x}; u)\langle J(\bar{x}), o(\bar{x}; u)\rangle|$$
$$= |\langle J(a(\bar{x}; u)\bar{x}), a(\bar{x}; u)\bar{x}\rangle|$$
$$= \|J(a(\bar{x}; u)\bar{x})\|_*^2$$
$$= \|a(\bar{x}; u)\bar{x}\|^2.$$

This proves (iv). Then, by part (iv), we prove part (v).

$$\|o(\bar{x}; u)\| = \|u - a(\bar{x}; u)\bar{x}\| \leq \|u\| + \|a(\bar{x}; u)\bar{x}\| \leq 2\|u\|. \qquad \square$$

Both $a(\bar{x}; \cdot)$ and $o(\bar{x}; \cdot)$ depend on $\bar{x}$. However, for the sake of simplicity, $a(\bar{x}; \cdot)$ and $o(\bar{x}; \cdot)$ are abbreviated as $a(\cdot)$ and $o(\cdot)$, respectively.

### 2.5. Strong smoothness of uniformly convex and uniformly smooth Banach spaces

In subsection 2.1, among the list of the properties of uniformly convex and uniformly smooth Banach space $X$, the property (xi) states that the norm of any uniformly convex and uniformly smooth Banach space is Gâteaux differentiable on the whole considered space. This property helps us to study the Gâteaux directional differentiability of mappings in uniformly convex and uniformly smooth Banach spaces (see [11]). The following results describe the strong smoothness of the norms of uniformly convex and uniformly smooth Banach spaces. This property will play important role to study the Gâteaux directional differentiability of the metric projection operator.

**Proposition 2.12.** *Let* $(X, \|\cdot\|)$ *be a real uniformly convex and uniformly smooth Banach space. Then, we have*

$$\lim_{\substack{v \to \theta \\ v \in O(\bar{x})}} \frac{\|\bar{x}+v\| - \|\bar{x}\|}{\|v\|} = 0, \text{ for each } \bar{x} \in X\setminus\{\theta\}. \tag{2.6}$$

*Proof.* For an arbitrarily fixed $\bar{x} \in X\setminus\{\theta\}$, by the uniformly convergence of (2.2) in Lemma 2.2, for any $\varepsilon > 0$, there is $\delta > 0$, such that, if $0 < t < \delta$, then

$$\left|\frac{\|x+ty\| - \|x\|}{t} - \langle J(x), y\rangle\right| < \varepsilon, \text{ for all } x, y \in \mathbb{S}.$$

This implies the following estimation. For any $v \in X\setminus\{\theta\}$ with $v \in O(\bar{x})$, if $\|v\| < \delta\|\bar{x}\|$, then, we have

$$\left|\frac{\|\bar{x}+v\|-\|\bar{x}\|}{\|v\|}\right|$$

$$= \left|\frac{\|\bar{x}\|\left(\left\|\frac{\bar{x}}{\|\bar{x}\|}+\frac{1}{\|\bar{x}\|}v\right\|-\left\|\frac{\bar{x}}{\|\bar{x}\|}\right\|\right)}{\|v\|}\right|$$

$$= \left|\frac{\left\|\frac{\bar{x}}{\|\bar{x}\|}+\frac{\|v\|}{\|\bar{x}\|}\frac{v}{\|v\|}\right\|-\left\|\frac{\bar{x}}{\|\bar{x}\|}\right\|}{\frac{\|v\|}{\|\bar{x}\|}}\right|$$

$$= \left|\frac{\left\|\frac{\bar{x}}{\|\bar{x}\|}+\frac{\|v\|}{\|\bar{x}\|}\frac{v}{\|v\|}\right\|-\left\|\frac{\bar{x}}{\|\bar{x}\|}\right\|}{\frac{\|v\|}{\|\bar{x}\|}} - \langle J\left(\frac{\bar{x}}{\|\bar{x}\|}\right), \frac{v}{\|v\|}\rangle\right|$$

$$< \varepsilon, \text{ for } \frac{\|v\|}{\|\bar{x}\|} < \delta \text{ with } v \in O(\bar{x}).$$

Here, since $v \in O(\bar{x})$, it implies that $\langle J\left(\frac{\bar{x}}{\|\bar{x}\|}\right), \frac{v}{\|v\|}\rangle = \frac{1}{\|\bar{x}\|\|v\|}\langle J(\bar{x}), v\rangle = 0$. □

## 3. The Fréchet differentiability of the metric projection onto closed balls

Theorem 3.3 in [13] proves the strict Fréchet differentiability of the metric projection onto closed balls in Hilbert spaces. However, the following theorem proves the Fréchet differentiability of the metric projection onto closed balls in uniformly convex and uniformly smooth Banach spaces. Before we prove next theorem, we review some related notations.

For any $x \in r\mathbb{S}$, two subsets $x_r^\uparrow$ and $x_r^\downarrow$ of $X\setminus\{\theta\}$ are defined by

(a) $x_r^\uparrow = \{v \in X\setminus\{\theta\}$: there is $\delta > 0$ such that $\|x + tv\| > r$, for all $t \in (0, \delta)\}$;
(b) $x_r^\downarrow = \{v \in X\setminus\{\theta\}$: there is $\delta > 0$ such that $\|x + tv\| \leq r$, for all $t \in (0, \delta)\}$.

**Theorem 3.1.** *Let X be a uniformly convex and uniformly smooth Banach space. For any $r > 0$, the metric projection $P_{r\mathbb{B}}: X \to r\mathbb{B}$ has the following differentiability properties.*

(i) $P_{r\mathbb{B}}$ *is strictly Fréchet differentiable on $r\mathbb{B}^o$ satisfying*

$$\nabla P_{r\mathbb{B}}(\bar{x}) = I_X, \text{ for every } \bar{x} \in r\mathbb{B}^o.$$

*That is,*

$$\bar{x} \in r\mathbb{B}^o \implies \nabla P_{r\mathbb{B}}(\bar{x})(x) = x, \text{ for every } x \in X;$$

(ii)    $P_{r\mathbb{B}}$ is Fréchet differentiable at every $\bar{x} \in X\setminus r\mathbb{B}$. The Fréchet derivative at $\bar{x}$ satisfies

$$\nabla P_{r\mathbb{B}}(\bar{x})(x) = \frac{r}{\|\bar{x}\|} o(\bar{x};\ x) = \frac{r}{\|\bar{x}\|}\left(x - \frac{\langle J(\bar{x}), x\rangle}{\|\bar{x}\|^2}\bar{x}\right),\ \text{for every } x \in X.$$

In particular, we have

(a) $\nabla P_{r\mathbb{B}}(\bar{x})(x) = \frac{r}{\|\bar{x}\|} x$, if $x \perp \bar{x}$, for $x \in X$;

(b) $\nabla P_{r\mathbb{B}}(\bar{x})(\bar{x}) = \theta$.

(iii)    For the subset $r\mathbb{S}$, we have

(I)    $P_{r\mathbb{B}}$ is Gâteaux directional differentiable on $r\mathbb{S}$ satisfying that, for every point $\bar{x} \in r\mathbb{S}$, the following representations are satisfied

(a)    $P'_{r\mathbb{B}}(\bar{x})(w) = w - \frac{1}{r^2}\langle J(\bar{x}), w\rangle \bar{x},\ \text{if } w \in \bar{x}_r^\uparrow$;

(b)    $P'_{r\mathbb{B}}(\bar{x})(\bar{x}) = \theta$;

(c)    $P'_{r\mathbb{B}}(\bar{x})(w) = w,\ \text{if } w \in \bar{x}_r^\downarrow$.

(II)    $P_{r\mathbb{B}}$ is not Fréchet differentiable at any point $\bar{x} \in r\mathbb{S}$. That is,

$$\nabla P_{r\mathbb{B}}(\bar{x})\ \text{does not exist, for any } \bar{x} \in r\mathbb{S}.$$

*Proof.* Proof of (i). Part (i) follows from Proposition 2.10. For any given $\bar{x} \in r\mathbb{B}^o$, there is $q > 0$ such that $\mathbb{B}(\bar{x}, q) \subseteq r\mathbb{B}^o$. By (2.4), we calculate

$$\lim_{u \to \bar{x}, v \to \bar{x}} \frac{P_{r\mathbb{B}}(u) - P_{r\mathbb{B}}(v) - I_X(u-v)}{\|u-v\|}$$

$$= \lim_{\substack{u \to \bar{x}, v \to \bar{x} \\ u,v \in \mathbb{B}(\bar{x},q)}} \frac{P_{r\mathbb{B}}(u) - P_{r\mathbb{B}}(v) - (u-v)}{\|u-v\|}$$

$$= \lim_{\substack{u \to \bar{x}, v \to \bar{x} \\ u,v \in \mathbb{B}(\bar{x},q)}} \frac{u - v - (u-v)}{\|u-v\|}$$

$$= \theta.$$

Hence, $P_{r\mathbb{B}}$ is strictly Fréchet differentiable at $\bar{x}$, with $\nabla P_{r\mathbb{B}}(\bar{x}) = I_X$, for $\bar{x} \in r\mathbb{B}^o$.

Proof of (ii). Let $\bar{x} \in X\setminus r\mathbb{B}$ be arbitrarily given with $\|\bar{x}\| > r$. By the definition of $o(\cdot)$, we may actually write

$$\nabla P_{r\mathbb{B}}(\bar{x})(x) = \frac{r}{\|\bar{x}\|}\left(x - \frac{\langle J(\bar{x}), x\rangle}{\|\bar{x}\|^2}\bar{x}\right) = \frac{r}{\|\bar{x}\|} o(x), \text{ for every } x \in X.$$

To prove part (ii) of this theorem, by $o(\bar{x}) = \theta$, we only need to verify that the above formula satisfies the following equations.

$$\theta = \lim_{u \to \bar{x}} \frac{P_{r\mathbb{B}}(u) - P_{r\mathbb{B}}(\bar{x}) - \frac{r}{\|\bar{x}\|} o(u - \bar{x})}{\|u - \bar{x}\|}$$

$$= \lim_{u \to \bar{x}} \frac{P_{r\mathbb{B}}(a(u)\bar{x} + o(u)) - P_{r\mathbb{B}}(\bar{x}) - \frac{r}{\|\bar{x}\|} o(u - \bar{x})}{\|u - \bar{x}\|}$$

$$= \lim_{u \to \bar{x}} \frac{P_{r\mathbb{B}}(a(u)\bar{x} + o(u)) - P_{r\mathbb{B}}(\bar{x}) - \frac{r}{\|\bar{x}\|} o(u) + \frac{r}{\|\bar{x}\|} o(\bar{x})}{\|u - \bar{x}\|}$$

$$= \lim_{u \to \bar{x}} \frac{P_{r\mathbb{B}}(a(u)\bar{x} + o(u)) - P_{r\mathbb{B}}(\bar{x}) - \frac{r}{\|\bar{x}\|} o(u)}{\|u - \bar{x}\|}$$

Since $\bar{x} \notin r\mathbb{B}$ with $\|\bar{x}\| > r$, there is $p > 0$ such that

$$\mathbb{B}(\bar{x}, p) \cap r\mathbb{B} = \emptyset.$$

By (2.4), this implies

$$\lim_{u \to \bar{x}} \frac{P_{r\mathbb{B}}(a(u)\bar{x} + o(u)) - P_{r\mathbb{B}}(\bar{x}) - \frac{r}{\|\bar{x}\|} o(u)}{\|u - \bar{x}\|}$$

$$= \lim_{\substack{u \to \bar{x} \\ u \in \mathbb{B}(\bar{x}, p)}} \frac{P_{r\mathbb{B}}(a(u)\bar{x} + o(u)) - P_{r\mathbb{B}}(\bar{x}) - \frac{r}{\|\bar{x}\|} o(u)}{\|u - \bar{x}\|}$$

$$= \lim_{\substack{u \to \bar{x} \\ u \in \mathbb{B}(\bar{x}, p)}} \frac{\frac{r}{\|a(u)\bar{x} + o(u)\|}(a(u)\bar{x} + o(u)) - \frac{r}{\|\bar{x}\|}\bar{x} - \frac{r}{\|\bar{x}\|} o(u)}{\|u - \bar{x}\|}$$

$$= r \lim_{\substack{u \to \bar{x} \\ u \in \mathbb{B}(\bar{x}, p)}} \frac{\frac{a(u)\bar{x} + o(u)}{\|a(u)\bar{x} + o(u)\|} - \frac{\bar{x}}{\|\bar{x}\|} - \frac{o(u)}{\|\bar{x}\|}}{\|u - \bar{x}\|}$$

$$= r \lim_{\substack{u \to \bar{x} \\ u \in \mathbb{B}(\bar{x}, p)}} \frac{\frac{a(u)\bar{x}}{\|a(u)\bar{x} + o(u)\|} + \frac{o(u)}{\|a(u)\bar{x} + o(u)\|} - \frac{o(u)}{\|\bar{x}\|} - \frac{\bar{x}}{\|\bar{x}\|}}{\|u - \bar{x}\|}$$

$$= r \lim_{\substack{u \to \bar{x} \\ u \in \mathbb{B}(\bar{x}, p)}} \frac{\frac{a(u)\bar{x}}{\|a(u)\bar{x} + o(u)\|} - \frac{\bar{x}}{\|\bar{x}\|} + \frac{o(u)}{\|a(u)\bar{x} + o(u)\|} - \frac{o(u)}{\|\bar{x}\|}}{\|u - \bar{x}\|}$$

$$= r \lim_{\substack{u \to \bar{x} \\ u \in \mathbb{B}(\bar{x}, p)}} \left( \frac{\frac{a(u)\bar{x}}{\|a(u)\bar{x} + o(u)\|} - \frac{\bar{x}}{\|\bar{x}\|}}{\|u - \bar{x}\|} + \frac{\frac{o(u)}{\|a(u)\bar{x} + o(u)\|} - \frac{o(u)}{\|\bar{x}\|}}{\|u - \bar{x}\|} \right). \tag{3.3}$$

By (iii) in Lemma 2.9, we have

$$u \to \bar{x} \quad \Longleftrightarrow \quad a(u) \to 1 \text{ and } o(u) \to \theta. \tag{3.4}$$

At first, we estimate the first part in the limit (3.3).

$$\left\| \frac{\frac{a(u)\bar{x}}{\|a(u)\bar{x}+o(u)\|} - \frac{\bar{x}}{\|\bar{x}\|}}{\|u-\bar{x}\|} \right\|$$

$$= \frac{\left\| \frac{a(u)\bar{x}}{\|a(u)\bar{x}+o(u)\|} - \frac{\bar{x}}{\|\bar{x}\|} \right\|}{\|u-\bar{x}\|}$$

$$= \frac{\|\bar{x}\| \left| \frac{a(u)}{\|a(u)\bar{x}+o(u)\|} - \frac{1}{\|\bar{x}\|} \right|}{\|u-\bar{x}\|}$$

$$= \frac{\|\bar{x}\| \left| \frac{a(u)\|\bar{x}\| - \|a(u)\bar{x}+o(u)\|}{\|a(u)\bar{x}+o(u)\|\|\bar{x}\|} \right|}{\|u-\bar{x}\|}$$

$$= \frac{\|\bar{x}\| \left| \frac{(a(u)\|\bar{x}\|)^2 - \|a(u)\bar{x}+o(u)\|^2}{\|a(u)\bar{x}+o(u)\|\,\|\bar{x}\|\,(a(u)\|\bar{x}\| + \|a(u)\bar{x}+o(u)\|)} \right|}{\|u-\bar{x}\|}$$

$$= \frac{\|\bar{x}\|}{\|a(u)\bar{x}+o(u)\|\,\|\bar{x}\|\,(a(u)\|\bar{x}\| + \|a(u)\bar{x}+o(u)\|)} \cdot \frac{\left| (a(u)\|\bar{x}\|)^2 - \|a(u)\bar{x}+o(u)\|^2 \right|}{\|u-\bar{x}\|}$$

$$= \frac{(a(u))^2 \|\bar{x}\|}{\|a(u)\bar{x}+o(u)\|\,\|\bar{x}\|\,(a(u)\|\bar{x}\| + \|a(u)\bar{x}+o(u)\|)} \cdot \frac{\left| \left\| \bar{x} + o\!\left(\frac{u}{a(u)}\right) \right\|^2 - \|\bar{x}\|^2 \right|}{\|u-\bar{x}\|}$$

By Property (v) in Lemma 2.9 and by the fact that $o(\bar{x}) = \theta$, we have

$$\|u - \bar{x}\| \geq \tfrac{1}{2}\|o(u-\bar{x})\| = \tfrac{1}{2}\|o(u)\|.$$

By (3.4) and Proposition 2.12, this implies

$$\frac{(a(u))^2\|\bar{x}\|}{\|a(u)\bar{x}+o(u)\|\,\|\bar{x}\|\,(a(u)\|\bar{x}\|+\|a(u)\bar{x}+o(u)\|)} \cdot \frac{\left|\left\|\bar{x}+o\!\left(\frac{u}{a(u)}\right)\right\|^2 - \|\bar{x}\|^2\right|}{\|u-\bar{x}\|}$$

$$\leq \frac{(a(u))^2\|\bar{x}\|}{\|a(u)\bar{x}+o(u)\|\,\|\bar{x}\|\,(a(u)\|\bar{x}\|+\|a(u)\bar{x}+o(u)\|)} \cdot \frac{\left|\left\|\bar{x}+o\!\left(\frac{u}{a(u)}\right)\right\|^2 - \|\bar{x}\|^2\right|}{\tfrac{1}{2}\|o(u)\|}$$

$$= \frac{2|a(u)|\|\bar{x}\|}{\|a(u)\bar{x}+o(u)\|\,\|\bar{x}\|\,(a(u)\|\bar{x}\|+\|a(u)\bar{x}+o(u)\|)} \cdot \frac{\left|\left\|\bar{x}+o\!\left(\frac{u}{a(u)}\right)\right\|^2 - \|\bar{x}\|^2\right|}{\left\|o\!\left(\frac{u}{a(u)}\right)\right\|}$$

$$\to \frac{1}{\|\bar{x}\|^2}\,0 = 0, \quad \text{as } u \to \bar{x}.$$

This proves the first part in (3.3)

$$\lim_{\substack{u\to\bar{x}\\ u\in\mathbb{B}(\bar{x},p)}} \frac{\frac{a(u)\bar{x}}{\|a(u)\bar{x}+o(u)\|} - \frac{\bar{x}}{\|\bar{x}\|}}{\|u-\bar{x}\|} = \theta. \tag{3.5}$$

Next, we similarly estimate the second part in (3.3).

$$\left\| \frac{\frac{o(u)}{\|a(u)\bar{x}+o(u)\|} - \frac{o(u)}{\|\bar{x}\|}}{\|u-\bar{x}\|} \right\|$$

$$= \frac{\|o(u)\| \left| \frac{1}{\|a(u)\bar{x}+o(u)\|} - \frac{1}{\|\bar{x}\|} \right|}{\|u-\bar{x}\|}$$

$$= \frac{\|o(u)\|}{\|a(u)\bar{x}+o(u)\|\|\bar{x}\|} \cdot \frac{|\|a(u)\bar{x}+o(u)\| - \|\bar{x}\||}{\|u-\bar{x}\|}$$

$$= \frac{|\|a(u)\bar{x}+o(u)\| - \|\bar{x}\||}{\|a(u)\bar{x}+o(u)\|\|\bar{x}\|} \cdot \frac{\|o(u)\|}{\|u-\bar{x}\|}$$

$$\leq \frac{|\|a(u)\bar{x}+o(u)\| - \|\bar{x}\||}{\|a(u)\bar{x}+o(u)\|\|\bar{x}\|} \cdot \frac{\|o(u)\|}{\frac{1}{2}\|o(u)\|}$$

$$= \frac{2|\|a(u)\bar{x}+o(u)\| - \|\bar{x}\||}{\|a(u)\bar{x}+o(u)\|\|\bar{x}\|}$$

$$\to \frac{0}{\|\bar{x}\|^2} = 0, \quad \text{as } u \to \bar{x}.$$

This implies

$$\lim_{\substack{u \to \bar{x} \\ u \in \mathbb{B}(\bar{x},p)}} \frac{\frac{o(u)}{\|a(u)\bar{x}+o(u)\|} - \frac{o(u)}{\|\bar{x}\|}}{\|u-\bar{x}\|} = \theta. \tag{3.6}$$

Combining (3.3), (3.5) and (3.6), we proved

$$\lim_{u \to \bar{x}} \frac{P_{r\mathbb{B}}(u) - P_{r\mathbb{B}}(\bar{x}) - \frac{r}{\|\bar{x}\|}o(u-\bar{x})}{\|u-\bar{x}\|} = \theta.$$

Proof of (II) in (iii). The proof of (II) is same with the proof of part (iii) in Theorem 3.3 in [13]. For an arbitrary given $\bar{x} \in r\mathbb{S}$, we assume, by the way of contradiction, that $P_{r\mathbb{B}}$ is Fréchet differentiable at $\bar{x}$. Then, there is a linear and continuous mapping $A(\bar{x}): X \to X$, such that

$$\lim_{u \to \bar{x}} \frac{P_{r\mathbb{B}}(u) - P_{r\mathbb{B}}(\bar{x}) - A(\bar{x})(u-\bar{x})}{\|u-\bar{x}\|} = \theta.$$

In particular, in the above limit, we take a directional line segment $u = (1+\delta)\bar{x}$, for $\delta \downarrow 0$. Since $A(\bar{x})$ is assumed to be linear and continuous, we have

$$\theta = \lim_{\delta \downarrow 0} \frac{P_{r\mathbb{B}}((1+\delta)\bar{x}) - P_{r\mathbb{B}}(\bar{x}) - A(\bar{x})((1+\delta)\bar{x}-\bar{x})}{\|(1+\delta)\bar{x}-\bar{x}\|}$$

$$= \lim_{\delta \downarrow 0} \frac{\frac{r}{\|(1+\delta)\bar{x}\|}(1+\delta)\bar{x} - \bar{x} - \delta A(\bar{x})(\bar{x})}{\delta\|\bar{x}\|}$$

$$= \lim_{\delta \downarrow 0} \frac{\frac{r}{(1+\delta)\|\bar{x}\|}(1+\delta)\bar{x} - \bar{x} - \delta A(\bar{x})(\bar{x})}{\delta\|\bar{x}\|}$$

$$= \lim_{\delta \downarrow 0} \frac{\bar{x}-\bar{x}}{\delta \|\bar{x}\|} - \frac{A(\bar{x})(\bar{x})}{\|\bar{x}\|}$$

$$= \theta - \frac{A(\bar{x})(\bar{x})}{r}$$

$$= -\frac{A(\bar{x})(\bar{x})}{r}.$$

This implies

$$A(\bar{x})(\bar{x}) = \theta. \tag{3.7}$$

Next, we take an opposite directional line segment $v = (1 - \delta)\bar{x}$, for $\delta \downarrow 0$ with $0 < \delta < 1$. Since $\|\bar{x}\| = r$, it follows that $(1 - \delta)\bar{x} \in r\mathbb{B}$, for any $\delta$ with $0 < \delta < 1$. By (2.4) and by the assumed linearity of $A(\bar{x})$, we have.

$$\theta = \lim_{v \to \bar{x}} \frac{P_{r\mathbb{B}}(v) - P_{r\mathbb{B}}(\bar{x}) - A(\bar{x})(v-\bar{x})}{\|x-\bar{x}\|}$$

$$= \lim_{\delta \downarrow 0, \delta < 1} \frac{(1-\delta)\bar{x} - \bar{x} - A(\bar{x})((1-\delta)\bar{x} - \bar{x})}{\|(1-\delta)\bar{x} - \bar{x}\|}$$

$$= \lim_{\delta \downarrow 0, \delta < 1} \frac{-\delta\bar{x} + \delta A(\bar{x})(\bar{x})}{\delta\|\bar{x}\|}$$

$$= \frac{-\bar{x} + A(\bar{x})(\bar{x})}{r}.$$

This implies

$$A(\bar{x})(\bar{x}) = \bar{x}.$$

By $\bar{x} \neq \theta$, this contradicts to (3.7), which proves that $P_{r\mathbb{B}}$ is not Fréchet differentiable at any point $\bar{x} \in r\mathbb{S}$. Part (II) in (iii) of this theorem is proved.

Proof of part (I) in (iii). By part (iii) in Theorem 5.2 in [11], we have that $P_{r\mathbb{B}}$ is Gâteaux directional differentiable on $r\mathbb{S}$. For every given point $\bar{x} \in r\mathbb{S}$, taking $c = \theta$ in (a) of part (iii) in Theorem 5.2 in [11], we have

$$P'_{r\mathbb{B}}(\bar{x})(w) = w - \frac{\|w\|}{r} \psi(\frac{\bar{x}}{r}, \frac{w}{\|w\|})\bar{x}, \quad \text{if } w \in \bar{x}_r^{\uparrow}. \tag{3.8}$$

By (2.2) in Lemma 2.2, we have

$$\psi(\frac{\bar{x}}{r}, \frac{w}{\|w\|}) = \langle J\left(\frac{\bar{x}}{r}\right), \frac{w}{\|w\|} \rangle = \frac{1}{r\|w\|} \langle J(\bar{x}), w \rangle.$$

Substituting this into (3.8), we obtain

$$P'_{r\mathbb{B}}(\bar{x})(w) = w - \frac{1}{r^2} \langle J(\bar{x}), w \rangle \bar{x}, \quad \text{if } w \in \bar{x}_r^{\uparrow}.$$

This proves (a) of part (I) in (iii) of this theorem. Then, (b) of (I) follows immediately from (a) of

(I). (c) of (I) in part (iii) of this theorem follows from (b) of part (iii) in Theorem 5.2 in [11]. This proves part (I) in (iii) of this theorem.

Proof of part (II) of (iii). The proof of part (II) of (iii) in this theorem is same with the proof of part (II) of (iii) in Theorem 3.3 in [13], which is omitted here. □

**Example 3.2.** We consider the real uniformly convex and uniformly smooth Banach space $(l_p, \|\cdot\|_p)$ with dual space $(l_q, \|\cdot\|_q)$, where both $p$ and $q$ are positive satisfying $1 < p, q < \infty$ and $\frac{1}{p} + \frac{1}{q} = 1$. The normalized duality mapping $J: l_p \to l_q$ has the following representations, for any $z = (z_1, z_2, \ldots) \in l_p$ with $z \neq \theta$,

$$Jz = \left( \frac{|z_1|^{p-1}\text{sign}(z_1)}{\|z\|_p^{p-2}}, \frac{|z_2|^{p-1}\text{sign}(z_2)}{\|z\|_p^{p-2}}, \ldots \right)$$

$$= \left( \frac{|z_1|^{p-2}z_1}{\|z\|_p^{p-2}}, \frac{|z_2|^{p-2}z_2}{\|z\|_p^{p-2}}, \ldots \right).$$

For any given $r > 0$, by Theorem 3.1, $P_{r\mathbb{B}}: l_p \to r\mathbb{B}$ has the following differentiability properties.

(i) $P_{r\mathbb{B}}$ is strictly Fréchet differentiable on $r\mathbb{B}^\circ$ such that, for any $\bar{x} \in r\mathbb{B}^\circ$, we have

$$\nabla P_{r\mathbb{B}}(\bar{x})(x) = x, \text{ for every } x \in l_p;$$

(ii) $P_{r\mathbb{B}}$ is Fréchet differentiable at every $\bar{x} \in l_p \setminus r\mathbb{B}$ such that, for every $x \in l_p$, we have

$$\nabla P_{r\mathbb{B}}(\bar{x})(x) = \frac{r}{\|\bar{x}\|_p}\left( x - \frac{\langle J(\bar{x}), x \rangle}{\|\bar{x}\|_p^2} \bar{x} \right) = \frac{r}{\|\bar{x}\|_p}\left( x - \frac{\sum_{i=1}^\infty |\bar{x}_i|^{p-2} \bar{x}_i x_i}{\|\bar{x}\|_p^p} \bar{x} \right).$$

(iii) For every point $\bar{x} \in r\mathbb{S}$, $\nabla P_{r\mathbb{B}}(\bar{x})$ does not exist. However, $P_{r\mathbb{B}}$ is Gâteaux directional differentiable at $\bar{x} \in r\mathbb{S}$, such that

(a) $P'_{r\mathbb{B}}(\bar{x})(w) = w - \frac{1}{r^2} \langle J(\bar{x}), w \rangle \bar{x} = w - \frac{\sum_{i=1}^\infty |\bar{x}_i|^{p-2} \bar{x}_i w_i}{r^p} \bar{x}$, if $w \in \bar{x}_r^\uparrow$;

(b) $P'_{r\mathbb{B}}(\bar{x})(\bar{x}) = \theta$;

(c) $P'_{r\mathbb{B}}(\bar{x})(w) = w$, if $w \in \bar{x}_r^\downarrow$.

Next example can be considered as a special case of Example 3.2.

**Example 3.3.** Let $X = \mathbb{R}^3$, where $(\mathbb{R}^3, \|\cdot\|_3)$ is the uniformly convex and uniformly smooth Banach space $\mathbb{R}^3$ equipped with the 3-norm $\|\cdot\|_3$, which is used in Proposition 2.5. Recall that, for any point $z = (z_1, z_2, z_3) \in \mathbb{R}^3$,

$$\|z\|_3 = \sqrt[3]{|z_1|^3 + |z_2|^3 + |z_3|^3}.$$

The dual space of $(\mathbb{R}^3, \|\cdot\|_3)$ is $(\mathbb{R}^3, \|\cdot\|_{\frac{3}{2}})$. The normalized duality mapping $J: \mathbb{R}^3 \to \mathbb{R}^3$ has the following representations, for any $z = (z_1, z_2, z_3) \in \mathbb{R}^3$ with $z \neq \theta$,

$$Jz = \left(\frac{|z_1|^2 \operatorname{sign}(z_1)}{\|z\|_3}, \frac{|z_2|^2 \operatorname{sign}(z_2)}{\|z\|_3}, \frac{|z_3|^2 \operatorname{sign}(z_3)}{\|z\|_3}\right)$$

$$= \left(\frac{|z_1|z_1}{\|z\|_3}, \frac{|z_2|z_2}{\|z\|_3}, \frac{|z_3|z_3}{\|z\|_3}\right).$$

For any given $r > 0$, by Theorem 3.1, $P_{r\mathbb{B}}: \mathbb{R}^3 \to r\mathbb{B}$ has the following differentiability properties.

(i) $P_{r\mathbb{B}}$ is strictly Fréchet differentiable on $r\mathbb{B}^o$ such that, for any $\bar{x} \in r\mathbb{B}^o$, we have

$$\nabla P_{r\mathbb{B}}(\bar{x})(x) = x, \text{ for every } x \in \mathbb{R}^3;$$

(ii) $P_{r\mathbb{B}}$ is Fréchet differentiable at every $\bar{x} \in \mathbb{R}^3 \setminus r\mathbb{B}$ such that, for every $x \in \mathbb{R}^3$, we have

$$\nabla P_{r\mathbb{B}}(\bar{x})(x) = \frac{r}{\|\bar{x}\|_3}\left(x - \frac{\langle J(\bar{x}), x\rangle}{\|\bar{x}\|_3^2}\bar{x}\right) = \frac{r}{\|\bar{x}\|_3}\left(x - \frac{\sum_{i=1}^{3}|\bar{x}_i|\bar{x}_i x_i}{\|\bar{x}\|_3^3}\bar{x}\right).$$

(iii) For every point $\bar{x} \in r\mathbb{S}$, $\nabla P_{r\mathbb{B}}(\bar{x})$ does not exist. However, $P_{r\mathbb{B}}$ is Gâteaux directional differentiable at $\bar{x} \in r\mathbb{S}$, such that

(a) $P'_{r\mathbb{B}}(\bar{x})(w) = w - \frac{1}{r^2}\langle J(\bar{x}), w\rangle \bar{x} = w - \frac{\sum_{i=1}^{3}|\bar{x}_i|\bar{x}_i w_i}{r^3}\bar{x}, \quad \text{if } w \in \bar{x}_r^\uparrow;$

(b) $P'_{r\mathbb{B}}(\bar{x})(\bar{x}) = \theta;$

(c) $P'_{r\mathbb{B}}(\bar{x})(w) = w, \quad \text{if } w \in \bar{x}_r^\downarrow.$

## 4. The Fréchet differentiability of the metric projection onto closed and convex cylinders in real Banach space $l_p$

In this section, we focus on the real uniformly convex and uniformly smooth Banach spaces $(l_p, \|\cdot\|_p)$ and $(l_q, \|\cdot\|_q)$ satisfying $1 < p, q < \infty$ and $\frac{1}{p} + \frac{1}{q} = 1$. As mentioned in Example 3.2, both $l_p$ and $l_q$ have origin $\theta = (0, 0, \ldots)$. They are the dual spaces of each other. We define some closed and convex cylinders in $l_p$. Then, we investigate the Fréchet differentiability of the metric projection onto closed and convex cylinders in $l_p$.

Recall the representations of the normalized duality mapping $J: l_p \to l_q$ given in Example 3.2. For any point $x = (x_1, x_2, \ldots) \in l_p$ with $x \neq \theta$, we have

$$J(x) = \left(\frac{|x_1|^{p-1}\operatorname{sign}(x_1)}{\|x\|_p^{p-2}}, \frac{|x_2|^{p-1}\operatorname{sign}(x_2)}{\|x\|_p^{p-2}}, \ldots\right)$$

$$= \left(\frac{|x_1|^{p-2}x_1}{\|x\|_p^{p-2}}, \frac{|x_2|^{p-2}x_2}{\|x\|_p^{p-2}}, \ldots\right). \tag{4.1}$$

Similarly, to (4.1), the representations of the normalized duality mapping $J^*: l_q \to l_p$ is given, for any $y = (y_1, y_2, \ldots) \in l_q$ with $y \neq \theta$ by

$$J^*(y) = \left(\frac{|y_1|^{q-1}\operatorname{sign}(y_1)}{\|y\|_q^{q-2}}, \frac{|y_2|^{q-1}\operatorname{sign}(y_2)}{\|y\|_q^{q-2}}, \ldots\right)$$

$$= \left( \frac{|y_1|^{q-2} y_1}{\|y\|_q^{q-2}}, \frac{|y_2|^{q-2} y_2}{\|y\|_q^{q-2}}, \dots \right). \tag{4.2}$$

Let $\mathbb{N}$ denote the set of all positive integers. Let $M$ be a nonempty subset of $\mathbb{N}$ with complementary set $\bar{M} = \mathbb{N} \backslash M$. We define

$$l_p^M = \{x = (x_1, x_2, \dots) \in l_p : x_i = 0, \text{ for all } i \in \bar{M}\},$$

$$l_q^M = \{y = (y_1, y_2, \dots) \in l_q : y_i = 0, \text{ for all } i \in \bar{M}\}.$$

$l_p^M$ and $l_q^M$ are closed subspaces of $l_p$ and $l_q$, respectively. They are the duality spaces of each other. We define a mapping $(\cdot)_M : l_p \to l_p^M$, for $x = (x_1, x_2, \dots) \in l_p$ by

$$(x_M)_i = \begin{cases} x_i, & \text{for } i \in M, \\ 0, & \text{for } i \notin M, \end{cases} \text{ for } i \in \mathbb{N}.$$

Similarly, we define a mapping $(\cdot)_{\bar{M}} : l_p \to l_p^{\bar{M}}$, for $x = (x_1, x_2, \dots) \in l_p$ by

$$(x_{\bar{M}})_i = \begin{cases} x_i, & \text{for } i \in \bar{M}, \\ 0, & \text{for } i \notin \bar{M}, \end{cases} \text{ for } i \in \mathbb{N}.$$

Then, $l_p$ has the following decomposition

$$x = x_M + x_{\bar{M}}, \text{ for any } x \in l_p. \tag{4.3}$$

**Lemma 4.1**. *Let $M$ be a nonempty subset of $\mathbb{N}$. Then $J$ is the normalized duality mapping from $l_p^M$ to $l_q^M$. That is,*

$$J(x) \in l_q^M, \text{ for any } x \in l_p^M.$$

*Proof*. The proof of this lemma follows immediately from (4.1). □

Let $\mathbb{B}_M$ denote the unit closed ball in $l_p^M$. It follows that, for any $r > 0$, $r\mathbb{B}_M$ is the closed ball with radius $r$ and centered at the origin in $l_p^M$. Let $\mathbb{S}_M$ be the unit sphere in $l_p^M$. Then, $r\mathbb{S}_M$ is the sphere in $l_p^M$ with radius $r$ and center $\theta$. We define

$$\mathbb{C}_M = \{x \in l_p : x_M \in \mathbb{B}_M\}.$$

$\mathbb{C}_M$ is called the cylinder in $l_p$ with base $\mathbb{B}_M$. It is a closed and convex subset in $l_p$. For any $r > 0$, $r\mathbb{C}_M$ is called the cylinder in $l_p$ with base $r\mathbb{B}_M$, which is a closed and convex subset of $l_p$. More precisely speaking, we have

$$r\mathbb{C}_M = \{x \in l_p : x_M \in r\mathbb{B}_M\}.$$

This implies that, for any $x \in l_p$, we have

$$x \in r\mathbb{C}_M \quad \Leftrightarrow \quad \|x_M\|_p \leq r. \tag{4.4}$$

The boundary of $r\mathbb{C}_M$ is denoted by $\partial(r\mathbb{C}_M)$ satisfying

$$\partial(r\mathbb{C}_M) = \{x \in l_p : \|x_M\|_p = r\}.$$

**Lemma 4.2.** *For any $r > 0$, the metric projection $P_{r\mathbb{C}_M}: l_p \to r\mathbb{C}_M$ satisfies the following formula.*

$$P_{r\mathbb{C}_M}(x) = \begin{cases} x, & \text{for any } x \in r\mathbb{C}_M, \\ \dfrac{r}{\|x_M\|_p} x_M + x_{\bar{M}}, & \text{for any } x \in l_p \setminus r\mathbb{C}_M. \end{cases} \quad (4.5)$$

*Proof.* It is clear that $P_{r\mathbb{C}_M}(x) = x$, for any $x \in r\mathbb{C}_M$, which proves the first part of (4.5). Next, we prove the second part of (4.5). By (4.4), we have

$$\|x_M\|_p > r, \text{ for any } x \in l_p \setminus r\mathbb{C}_M.$$

For any $z \in r\mathbb{C}_M$, by Lemma 4.1 and by (4.4), we calculate

$$\left\langle J\left(x - \left(\frac{r}{\|x_M\|_p} x_M + x_{\bar{M}}\right)\right), \frac{r}{\|x_M\|_p} x_M + x_{\bar{M}} - z \right\rangle$$

$$= \left\langle J\left((x_M + x_{\bar{M}}) - \left(\frac{r}{\|x_M\|_p} x_M + x_{\bar{M}}\right)\right), \frac{r}{\|x_M\|_p} x_M + x_{\bar{M}} - (z_M + z_{\bar{M}}) \right\rangle$$

$$= \left\langle J\left(\left(1 - \frac{r}{\|x_M\|_p}\right) x_M\right), \frac{r}{\|x_M\|_p} x_M - z_M + x_{\bar{M}} - z_{\bar{M}} \right\rangle$$

$$= \left(1 - \frac{r}{\|x_M\|_p}\right) \left\langle J(x_M), \frac{r}{\|x_M\|_p} x_M - z_M \right\rangle + \left(1 - \frac{r}{\|x_M\|_p}\right) \left\langle J(x_M), x_{\bar{M}} - z_{\bar{M}} \right\rangle$$

$$= \left(1 - \frac{r}{\|x_M\|_p}\right) \left\langle J(x_M), \frac{r}{\|x_M\|_p} x_M - z_M \right\rangle$$

$$= \left(1 - \frac{r}{\|x_M\|_p}\right) \left(\frac{r}{\|x_M\|_p} \|x_M\|_p^2 - \langle J(x_M), z_M \rangle\right)$$

$$\geq \left(1 - \frac{r}{\|x_M\|_p}\right) \left(r \|x_M\|_p - \|J(x_M)\|_q \|z_M\|_p\right)$$

$$= \left(1 - \frac{r}{\|x_M\|_p}\right) \left(r \|x_M\|_p - \|x_M\|_p \|z_M\|_p\right)$$

$$= \left(\|x_M\|_p - r\right)\left(r - \|z_M\|_p\right)$$

$\geq 0$, for all $z \in r\mathbb{C}_M$.

By the basic variational principle of $P_{r\mathbb{C}_M}$, this implies

$$P_{r\mathbb{C}_M}(x) = \frac{r}{\|x_M\|_p} x_M + x_{\bar{M}}, \text{ for any } x \in l_p \setminus r\mathbb{C}_M. \qquad \square$$

For any $x \in l_p$ with $\|x_M\|_p = r$, two subsets $x_r^{\Uparrow}$ and $x_r^{\Downarrow}$ of $l_p$ are defined by

(a) $x_r^{\Uparrow} = \{v \in l_p: \text{there is } \delta > 0 \text{ such that } \|(x+tv)_M\|_p > r, \text{ for all } t \in (0, \delta)\}$;
(b) $x_r^{\Downarrow} = \{v \in l_p: \text{there is } \delta > 0 \text{ such that } \|(x+tv)_M\|_p \le r, \text{ for all } t \in (0, \delta)\}$.

**Theorem 4.3.** *For any $r > 0$, the metric projection $P_{r\mathbb{C}_M}: l_p \to r\mathbb{C}_M$ has the following differentiability properties.*

(i) $P_{r\mathbb{C}_M}$ *is strictly Fréchet differentiable on $(r\mathbb{C}_M)^\circ$ satisfying*

$$\nabla P_{r\mathbb{C}_M}(\bar{x}) = I_{l_p}, \text{ for every } \bar{x} \in (r\mathbb{C}_M)^\circ.$$

*That is,*

$$\bar{x} \in (r\mathbb{C}_M)^\circ \implies \nabla P_{r\mathbb{C}_M}(\bar{x})(u) = u, \text{ for every } u \in l_p.$$

(ii) $P_{r\mathbb{C}_M}$ *is Fréchet differentiable at every point $\bar{x} \in l_p \setminus r\mathbb{C}_M$ such that*

$$\nabla P_{r\mathbb{C}_M}(\bar{x})(x) = \frac{r}{\|\bar{x}_M\|_p}\left(u_M - \frac{\langle J(\bar{x}_M), u_M \rangle}{\|\bar{x}_M\|_p^2}\bar{x}_M\right) + u_{\bar{M}}, \text{ for every } u \in l_p. \quad (4.6)$$

*In particular,*

$$\nabla P_{r\mathbb{C}_M}(\bar{x})(\bar{x}) = x_{\bar{M}}, \text{ for every } \bar{x} \in l_p \setminus r\mathbb{C}_M. \quad (4.7)$$

(iii) *On $\partial(r\mathbb{C}_M)$ (that is, $\bar{x} \in l_p$ with $\|\bar{x}_M\|_p = r$), we have*

(I) $P_{r\mathbb{C}_M}$ *is Gâteaux differentiable at every point $\bar{x} \in \partial(r\mathbb{C}_M)$ such that*

(a) $P'_{r\mathbb{C}_M}(\bar{x})(u) = u - \frac{\langle J(\bar{x}_M), u_M \rangle}{r^2}\bar{x}_M$, *if* $u \in \bar{x}_r^{\Uparrow}$;

(b) $P'_{r\mathbb{C}_M}(\bar{x})(\bar{x}) = \bar{x}_{\bar{M}}$;

(c) $P'_{r\mathbb{C}_M}(\bar{x})(u) = u$, *if* $u \in \bar{x}_r^{\Downarrow}$.

(II) $P_{r\mathbb{C}_M}$ *is not Fréchet differentiable at any point $\bar{x}$, that is,*

$$\nabla P_{r\mathbb{C}_M}(\bar{x}) \text{ does not exist, for any } \bar{x} \in \partial(r\mathbb{C}_M).$$

*Proof.* Part (i) follows from Proposition 2.10 immediately.

Proof of (ii). For any given $\bar{x} \in l_p \setminus r\mathbb{C}_M$, we prove

$$\theta = \lim_{u \to \bar{x}} \frac{P_{r\mathbb{C}_M}(u) - P_{r\mathbb{C}_M}(\bar{x}) - \nabla P_{r\mathbb{C}_M}(\bar{x})(u - \bar{x})}{\|u - \bar{x}\|}$$

$$= \lim_{u \to \bar{x}} \frac{P_{r\mathbb{C}_M}(u) - P_{r\mathbb{C}_M}(\bar{x}) - \left( \frac{r}{\|\bar{x}_M\|_p} \left( (u-\bar{x})_M - \frac{\langle J(\bar{x}_M),\, (u-\bar{x})_M \rangle}{\|\bar{x}_M\|_p^2} \bar{x}_M \right) + (u-\bar{x})_{\overline{M}} \right)}{\|u - \bar{x}\|}$$

$$= \lim_{u \to \bar{x}} \frac{P_{r\mathbb{C}_M}(u) - P_{r\mathbb{C}_M}(\bar{x}) - \left( \frac{r}{\|\bar{x}_M\|_p} \left( u_M - \frac{\langle J(\bar{x}_M),\, u_M \rangle}{\|\bar{x}_M\|_p^2} \bar{x}_M \right) + u_{\overline{M}} - \bar{x}_{\overline{M}} \right)}{\|u - \bar{x}\|}. \tag{4.8}$$

For this given $\bar{x} \in l_p \backslash r\mathbb{C}_M$, for any $u \in l_p$, we write

$$a(u_M) = \frac{\langle J(\bar{x}_M), u_M \rangle}{\|\bar{x}_M\|_p^2} \quad \text{and} \quad o(u_M) = u_M - \frac{\langle J(\bar{x}_M), u_M \rangle}{\|\bar{x}_M\|_p^2} \bar{x}_M.$$

Notice that both $a(u_M)$ and $o(u_M)$ depend on $\bar{x}_M$. Then, for any $u \in l_p$, $u_M$ enjoys the following decomposition

$$u_M = a(u_M)\bar{x}_M + o(u_M), \text{ for any } u \in l_p. \tag{4.9}$$

Since

$$u_M \to \bar{x}_M, \text{ as } u \to \bar{x}.$$

By the continuity of the normalized duality mapping $J$, this implies

$$a(u_M) \to 1 \quad \text{and} \quad o(u_M) \to \theta, \quad \text{as } u \to \bar{x}. \tag{4.10}$$

For this given $\bar{x} \in l_p \backslash r\mathbb{C}_M$, $a(u_M)$ and $o(u)$ has the following semi-orthogonal property

$$\langle J(a(u_M)\bar{x}_M),\, o(u_M) \rangle = 0, \quad \text{for any } u \in l_p. \tag{4.11}$$

Proof of (4.11). For any $u \in l_p$, by $\langle J(\bar{x}_M), \bar{x}_M \rangle = \|\bar{x}_M\|_p^2$, we calculate

$$\langle J(a(u_M)\bar{x}_M),\, o(u) \rangle$$

$$= \langle J(a(u_M)\bar{x}_M),\, u_M - \frac{\langle J(\bar{x}_M), u_M \rangle}{\|\bar{x}_M\|_p^2} \bar{x}_M \rangle$$

$$= a(u_M) \langle J(\bar{x}_M),\, u_M - \frac{\langle J(\bar{x}_M), u_M \rangle}{\|\bar{x}_M\|_p^2} \bar{x}_M \rangle$$

$$= a(u_M) \langle J(\bar{x}_M),\, u_M - \frac{\langle J(\bar{x}_M), u_M \rangle}{\|\bar{x}_M\|_p^2} \bar{x}_M \rangle$$

$$= a(u_M)(\langle J(\bar{x}_M), u_M \rangle - \langle J(\bar{x}_M), u_M \rangle)$$

$$= 0.$$

This proves (4.11). Now, we prove (4.8). Since $\bar{x} \in l_p \backslash r\mathbb{C}_M$ with $\|\bar{x}_M\|_p > r$, there is $b > 0$ such that, for any $u \in l_p$, we have

$$\|u - \bar{x}\|_p < b \implies \|u_M\|_p > r, \text{ which means } u \in l_p \setminus r\mathbb{C}_M.$$

By (4.5) and (4.9), this implies

$$\lim_{u \to \bar{x}} \frac{P_{r\mathbb{C}_M}(u) - P_{r\mathbb{C}_M}(\bar{x}) - \left(\frac{r}{\|\bar{x}_M\|_p}\left(u_M - \frac{\langle J(\bar{x}_M), u_M \rangle}{\|\bar{x}_M\|_p^2} \bar{x}_M\right) + u_{\overline{M}} - \bar{x}_{\overline{M}}\right)}{\|u - \bar{x}\|}$$

$$= \lim_{\substack{u \to \bar{x} \\ \|u - \bar{x}\|_p < b}} \frac{P_{r\mathbb{C}_M}(u) - P_{r\mathbb{C}_M}(\bar{x}) - \left(\frac{r}{\|\bar{x}_M\|_p}\left(u_M - \frac{\langle J(\bar{x}_M), u_M \rangle}{\|\bar{x}_M\|_p^2} \bar{x}_M\right) + u_{\overline{M}} - \bar{x}_{\overline{M}}\right)}{\|u - \bar{x}\|}$$

$$= \lim_{\substack{u \to \bar{x} \\ \|u - \bar{x}\|_p < b}} \frac{\left(\frac{r}{\|u_M\|_p} u_M + u_{\overline{M}}\right) - \left(\frac{r}{\|\bar{x}_M\|_p} \bar{x}_M + \bar{x}_{\overline{M}}\right) - \left(\frac{r}{\|\bar{x}_M\|_p}\left(u_M - \frac{\langle J(\bar{x}_M), u_M \rangle}{\|\bar{x}_M\|_p^2} \bar{x}_M\right) + u_{\overline{M}} - \bar{x}_{\overline{M}}\right)}{\|u - \bar{x}\|_p}$$

$$= \lim_{\substack{u \to \bar{x} \\ \|u - \bar{x}\|_p < b}} \frac{\frac{r}{\|u_M\|_p} u_M - \frac{r}{\|\bar{x}_M\|_p} \bar{x}_M - \frac{r}{\|\bar{x}_M\|_p}\left(u_M - \frac{\langle J(\bar{x}_M), u_M \rangle}{\|\bar{x}_M\|_p^2} \bar{x}_M\right)}{\|u - \bar{x}\|_p}$$

$$= \lim_{\substack{u \to \bar{x} \\ \|u - \bar{x}\|_p < b}} \frac{\frac{r}{\|a(u_M)\bar{x}_M + o(u_M)\|_p}(a(u_M)\bar{x}_M + o(u_M)) - \frac{r}{\|\bar{x}_M\|_p} \bar{x}_M - \frac{r}{\|\bar{x}_M\|_p} o(u_M)}{\|u - \bar{x}\|_p}$$

$$= r \lim_{\substack{u \to \bar{x} \\ \|u - \bar{x}\|_p < b}} \frac{\frac{a(u_M)\bar{x}_M}{\|a(u_M)\bar{x}_M + o(u_M)\|_p} - \frac{\bar{x}_M}{\|\bar{x}_M\|_p} + \frac{o(u_M)}{\|a(u_M)\bar{x}_M + o(u_M)\|_p} - \frac{o(u_M)}{\|\bar{x}_M\|_p}}{\|u - \bar{x}\|_p}$$

$$= r \lim_{\substack{u \to \bar{x} \\ \|u - \bar{x}\|_p < b}} \frac{\frac{a(u_M)\bar{x}_M}{\|a(u_M)\bar{x}_M + o(u_M)\|_p} - \frac{\bar{x}_M}{\|\bar{x}_M\|_p}}{\|u - \bar{x}\|_p} + r \lim_{\substack{u \to \bar{x} \\ \|u - \bar{x}\|_p < b}} \frac{\frac{o(u_M)}{\|a(u_M)\bar{x}_M + o(u_M)\|_p} - \frac{o(u_M)}{\|\bar{x}_M\|_p}}{\|u - \bar{x}\|_p}. \quad (4.12)$$

We estimate the first part of (4.12).

$$\lim_{\substack{u \to \bar{x} \\ \|u - \bar{x}\|_p < b}} \left\| \frac{\frac{a(u_M)\bar{x}_M}{\|a(u_M)\bar{x}_M + o(u_M)\|_p} - \frac{\bar{x}_M}{\|\bar{x}_M\|_p}}{\|u - \bar{x}\|_p} \right\|_p$$

$$= \lim_{\substack{u \to \bar{x} \\ \|u - \bar{x}\|_p < b}} \frac{\|\bar{x}_M\|_p}{\|a(u_M)\bar{x}_M + o(u_M)\|_p \|\bar{x}_M\|_p} \cdot \frac{\left|\|a(u_M)\bar{x}_M + o(u_M)\|_p - a(u_M)\|\bar{x}_M\|_p\right|}{\|u - \bar{x}\|_p}$$

$$= \frac{1}{\|\bar{x}_M\|_p} \lim_{\substack{u \to \bar{x} \\ \|u - \bar{x}\|_p < b}} \frac{\left|\|a(u_M)\bar{x}_M + o(u_M)\|_p - a(u_M)\|\bar{x}_M\|_p\right|}{\|u - \bar{x}\|_p}. \quad (4.13)$$

By (4.10), since $a(u_M) \to 1$, as $u \to \bar{x}$, we can suppose $a(u_M) > 0$ in the limit (4.13). This implies that (4.13) becomes

$$\frac{1}{\|\bar{x}_M\|_p} \lim_{\substack{u \to \bar{x} \\ \|u-\bar{x}\|_p < b}} \frac{\big|\|a(u_M)\bar{x}_M + o(u_M)\|_p - a(u_M)\|\bar{x}_M\|_p\big|}{\|u-\bar{x}\|_p}$$

$$= \frac{1}{\|\bar{x}_M\|_p} \lim_{\substack{u \to \bar{x} \\ \|u-\bar{x}\|_p < b}} \frac{\big|\|a(u_M)\bar{x}_M + o(u_M)\|_p - \|a(u_M)\bar{x}_M\|_p\big|}{\|u-\bar{x}\|_p}. \tag{4.14}$$

We find the connections between $\|o(u_M)\|_p$ and $\|u-\bar{x}\|_p$. By definition, we have

$$\|o(u_M)\|_p = \left\| u_M - \frac{\langle J(\bar{x}_M), u_M \rangle}{\|\bar{x}_M\|_p^2} \bar{x}_M \right\|_p$$

$$= \left\| \frac{\|\bar{x}_M\|_p^2}{\|\bar{x}_M\|_p^2} u_M - \frac{\langle J(\bar{x}_M), u_M \rangle}{\|\bar{x}_M\|_p^2} \bar{x}_M \right\|_p$$

$$= \frac{\|\langle J(\bar{x}_M), \bar{x}_M \rangle u_M - \langle J(\bar{x}_M), u_M \rangle \bar{x}_M\|_p}{\|\bar{x}_M\|_p^2}$$

$$= \frac{\|\langle J(\bar{x}_M), \bar{x}_M - u_M \rangle u_M + \langle J(\bar{x}_M), u_M \rangle (u_M - \bar{x}_M)\|_p}{\|\bar{x}_M\|_p^2}$$

$$\leq \frac{\|J(\bar{x}_M)\|_q \|\bar{x}_M - u_M\|_p \|u_M\|_p + \|J(\bar{x}_M)\|_q \|u_M\|_p \|\bar{x}_M - u_M\|_p}{\|\bar{x}_M\|_p^2}$$

$$= \frac{2\|\bar{x}_M\|_p \|u_M\|_p \|\bar{x}_M - u_M\|_p}{\|\bar{x}_M\|_p^2}$$

$$= \frac{2\|u_M\|_p}{\|\bar{x}_M\|_p} \|\bar{x}_M - u_M\|_p. \tag{4.15}$$

By definition, we calculate

$$\|u - \bar{x}\|_p^p$$

$$= \|u_M + u_{\bar{M}} - \bar{x}_M - \bar{x}_{\bar{M}}\|_p^p$$

$$= \|u_M - \bar{x}_M\|_p^p + \|u_{\bar{M}} - \bar{x}_{\bar{M}}\|_p^p$$

$$\geq \|u_M - \bar{x}_M\|_p^p.$$

This implies

$$\|u - \bar{x}\|_p \geq \|\bar{x}_M - u_M\|_p. \tag{4.16}$$

Substituting (4.16) into (4.15), we obtain

$$\|o(u_M)\|_p \leq \tfrac{2\|u_M\|_p}{\|\bar{x}_M\|_p} \|\bar{x}_M - u_M\|_p \leq \tfrac{2\|u_M\|_p}{\|\bar{x}_M\|_p} \|u - \bar{x}\|_p.$$

This implies

$$\|u - \bar{x}\|_p \geq \tfrac{\|\bar{x}_M\|_p}{2\|u_M\|_p} \|o(u_M)\|_p. \tag{4.17}$$

Substituting (4.17) into (4.13) and (4.14), and by Lemma 2.2 and by (4.11), we have

$$\lim_{\substack{u \to \bar{x} \\ \|u-\bar{x}\|_p < b}} \left\| \frac{\frac{a(u_M)\bar{x}_M}{\|a(u_M)\bar{x}_M + o(u_M)\|_p} - \frac{\bar{x}_M}{\|\bar{x}_M\|_p}}{\|u-\bar{x}\|_p} \right\|_p$$

$$= \tfrac{1}{\|\bar{x}_M\|_p} \lim_{\substack{u \to \bar{x} \\ \|u-\bar{x}\|_p < b}} \frac{\big|\|a(u_M)\bar{x}_M + o(u_M)\|_p - \|a(u_M)\bar{x}_M\|_p\big|}{\|u-\bar{x}\|_p}$$

$$\leq \tfrac{1}{\|\bar{x}_M\|_p} \lim_{\substack{u \to \bar{x} \\ \|u-\bar{x}\|_p < b}} \frac{\big|\|a(u_M)\bar{x}_M + o(u_M)\|_p - \|a(u_M)\bar{x}_M\|_p\big|}{\frac{\|\bar{x}_M\|_p}{2\|u_M\|_p} \|o(u_M)\|_p}$$

$$\leq \tfrac{2}{\|\bar{x}_M\|_p} \lim_{\substack{u \to \bar{x} \\ \|u-\bar{x}\|_p < b}} \frac{\big|\|a(u_M)\bar{x}_M + o(u_M)\|_p - \|a(u_M)\bar{x}_M\|_p\big|}{\|o(u_M)\|_p}$$

$$= \tfrac{2}{\|\bar{x}_M\|_p} \frac{|\langle J(a(u_M)\bar{x}_M), o(u_M)\rangle|}{\|a(u_M)\bar{x}_M\|_p}$$

$$= 0, \quad \text{for any } u \in l_p. \tag{4.18}$$

Hence, we proved the first part of (4.12)

$$\lim_{\substack{u \to \bar{x} \\ \|u-\bar{x}\|_p < b}} \left\| \frac{\frac{a(u_M)\bar{x}_M}{\|a(u_M)\bar{x}_M + o(u_M)\|_p} - \frac{\bar{x}_M}{\|\bar{x}_M\|_p}}{\|u-\bar{x}\|_p} \right\|_p = 0. \tag{4.19}$$

Next, we estimate the second part of (4.12). By the proof of (4.18), we have

$$\lim_{\substack{u \to \bar{x} \\ \|u-\bar{x}\|_p < b}} \left\| \frac{\frac{o(u_M)}{\|a(u_M)\bar{x}_M + o(u_M)\|_p} - \frac{o(u_M)}{\|\bar{x}_M\|_p}}{\|u-\bar{x}\|_p} \right\|_p$$

$$= \lim_{\substack{u \to \bar{x} \\ \|u-\bar{x}\|_p < b}} \frac{\|o(u_M)\|_p}{\|a(u_M)\bar{x}_M + o(u_M)\|_p \|\bar{x}_M\|_p} \frac{\big|\|a(u_M)\bar{x}_M + o(u_M)\|_p - \|\bar{x}_M\|_p\big|}{\|u-\bar{x}\|_p}$$

$$= \lim_{\substack{u \to \bar{x} \\ \|u-\bar{x}\|_p < b}} \frac{\|o(u_M)\|_p}{\|a(u_M)\bar{x}_M + o(u_M)\|_p \|\bar{x}_M\|_p} \frac{\big|\|a(u_M)\bar{x}_M + o(u_M)\|_p - \|a(u_M)\bar{x}_M\|_p + \|a(u_M)\bar{x}_M\|_p - \|\bar{x}_M\|_p\big|}{\|u-\bar{x}\|_p}$$

$$\leq \lim_{\substack{u \to \bar{x} \\ \|u-\bar{x}\|_p < b}} \frac{\|o(u_M)\|_p}{\|a(u_M)\bar{x}_M + o(u_M)\|_p \|\bar{x}_M\|_p} \frac{\big|\|a(u_M)\bar{x}_M + o(u_M)\|_p - \|a(u_M)\bar{x}_M\|_p\big|}{\|u-\bar{x}\|_p}$$

$$+ \lim_{\substack{u \to \bar{x} \\ \|u-\bar{x}\|_p < b}} \frac{\|o(u_M)\|_p}{\|a(u_M)\bar{x}_M + o(u_M)\|_p \|\bar{x}_M\|_p} \frac{\big|\|a(u_M)\bar{x}_M\|_p - \|\bar{x}_M\|_p\big|}{\|u-\bar{x}\|_p}$$

$$= 0 + \lim_{\substack{u \to \bar{x} \\ \|u-\bar{x}\|_p < b}} \frac{\|o(u_M)\|_p}{\|a(u_M)\bar{x}_M + o(u_M)\|_p \|\bar{x}_M\|_p} \frac{|a(u_M) - 1|\|\bar{x}_M\|_p}{\|u-\bar{x}\|_p}. \tag{4.20}$$

By the definition of $a(u_M)$ and by (4.16), we have

$$|a(u_M) - 1|$$

$$= \left|\frac{\langle J(\bar{x}_M), u_M \rangle}{\|\bar{x}_M\|_p^2} - 1\right|$$

$$= \left|\frac{\langle J(\bar{x}_M), u_M \rangle - \langle J(\bar{x}_M), \bar{x}_M \rangle}{\|\bar{x}_M\|_p^2}\right|$$

$$= \frac{|\langle J(\bar{x}_M), u_M \rangle - \langle J(\bar{x}_M), \bar{x}_M \rangle|}{\|\bar{x}_M\|_p^2}$$

$$\leq \frac{\|J(\bar{x}_M)\|_q \|\bar{x}_M - u_M\|_p}{\|\bar{x}_M\|_p^2}$$

$$= \frac{1}{\|\bar{x}_M\|_p} \|\bar{x}_M - u_M\|_p$$

$$\leq \frac{1}{\|\bar{x}_M\|_p} \|u - \bar{x}\|_p. \tag{4.21}$$

Substituting (4.21) into (4.20) and by (4.10), we have

$$\lim_{\substack{u \to \bar{x} \\ \|u-\bar{x}\|_p < b}} \left\|\frac{\frac{o(u_M)}{\|a(u_M)\bar{x}_M + o(u_M)\|_p} - \frac{o(u_M)}{\|\bar{x}_M\|_p}}{\|u-\bar{x}\|_p}\right\|_p$$

$$\leq 0 + \lim_{\substack{u \to \bar{x} \\ \|u-\bar{x}\|_p < b}} \frac{\|o(u_M)\|_p}{\|a(u_M)\bar{x}_M + o(u_M)\|_p \|\bar{x}_M\|_p} \frac{|a(u_M)-1|\|\bar{x}_M\|_p}{\|u-\bar{x}\|_p}$$

$$\leq \lim_{\substack{u \to \bar{x} \\ \|u-\bar{x}\|_p < b}} \frac{\|o(u_M)\|_p}{\|a(u_M)\bar{x}_M + o(u_M)\|_p \|\bar{x}_M\|_p}$$

$$= 0. \tag{4.22}$$

Substituting (4.19) and (4.22) into (4.12), we obtain

$$\lim_{u \to \bar{x}} \frac{P_{r\mathbb{C}_M}(u) - P_{r\mathbb{C}_M}(\bar{x}) - \left( \frac{r}{\|\bar{x}_M\|_p} \left( u_M - \frac{\langle J(\bar{x}_M), u_M \rangle}{\|\bar{x}_M\|_p^2} \bar{x}_M \right) + u_{\overline{M}} - \bar{x}_{\overline{M}} \right)}{\|u - \bar{x}\|} = \theta.$$

This proves part (ii) of this Theorem.

Next, we prove part (a) of (I) in (iii). Let $\bar{x} \in l_p$ with $\|\bar{x}_M\|_p = r$. For any $u \in \bar{x}_r^{\Uparrow}$, there is $\delta > 0$ such that $\|(x + tv)_M\|_p > r$, for all $t \in (0, \delta)$. Then, by Lemma 4.2 and Lemma 2.2, we have

$$P'_{r\mathbb{C}_M}(\bar{x})(u) = \lim_{t \downarrow 0} \frac{P_{r\mathbb{C}_M}(\bar{x}+tu) - P_{r\mathbb{C}_M}(\bar{x})}{t}$$

$$= \lim_{\substack{t \downarrow 0 \\ t < \delta}} \frac{P_{r\mathbb{C}_M}(\bar{x}+tu) - P_{r\mathbb{C}_M}(\bar{x})}{t}$$

$$= \lim_{\substack{t \downarrow 0 \\ t < \delta}} \frac{\frac{r}{\|(\bar{x}+tu)_M\|_p}(\bar{x}+tu)_M + (\bar{x}+tu)_{\overline{M}} - \bar{x}}{t}$$

$$= \lim_{\substack{t \downarrow 0 \\ t < \delta}} \frac{\frac{r}{\|(\bar{x}+tu)_M\|_p}(\bar{x}+tu)_M + (\bar{x}+tu)_{\overline{M}} - (\bar{x}_M + \bar{x}_{\overline{M}})}{t}$$

$$= \lim_{\substack{t \downarrow 0 \\ t < \delta}} \frac{\frac{r}{\|(\bar{x}+tu)_M\|_p}(\bar{x}_M + tu_M) + (\bar{x}_{\overline{M}} + tu_{\overline{M}}) - \left(\frac{r}{\|\bar{x}_M\|_p}\bar{x}_M + \bar{x}_{\overline{M}}\right)}{t}$$

$$= \lim_{\substack{t \downarrow 0 \\ t < \delta}} \frac{\frac{r}{\|(\bar{x}+tu)_M\|_p}\bar{x}_M - \frac{r}{\|\bar{x}_M\|_p}\bar{x}_M + \frac{r}{\|(\bar{x}+tu)_M\|_p}tu_M + tu_{\overline{M}}}{t}$$

$$= \lim_{\substack{t \downarrow 0 \\ t < \delta}} \frac{\frac{r}{\|(\bar{x}+tu)_M\|_p}\bar{x}_M - \frac{r}{\|\bar{x}_M\|_p}\bar{x}_M}{t} + \frac{r}{\|\bar{x}_M\|_p}u_M + u_{\overline{M}}$$

$$= r\bar{x}_M \lim_{\substack{t \downarrow 0 \\ t < \delta}} \frac{-1}{\|(\bar{x}+tu)_M\|_p \|\bar{x}_M\|_p} \frac{\|(\bar{x}+tu)_M\|_p - \|\bar{x}_M\|_p}{t} + \frac{r}{\|\bar{x}_M\|_p}u_M + u_{\overline{M}}$$

$$= r\bar{x}_M \lim_{\substack{t \downarrow 0 \\ t < \delta}} \frac{-1}{\|(\bar{x}+tu)_M\|_p \|\bar{x}_M\|_p} \frac{\|\bar{x}_M + tu_M\|_p - \|\bar{x}_M\|_p}{t} + u_M + u_{\overline{M}}$$

$$= r\bar{x}_M \frac{-1}{\|\bar{x}_M\|_p^2} \frac{\langle J(\bar{x}_M), u_M \rangle}{\|\bar{x}_M\|_p} + u_M + u_{\overline{M}}$$

$$= \bar{x}_M \frac{-\langle J(\bar{x}_M), u_M \rangle}{r^2} + u$$

$$= u - \frac{\langle J(\bar{x}_M), u_M \rangle}{r^2} \bar{x}_M.$$

This proves part (a) of (I) in (iii). Next, we prove (b). For any $\bar{x} \in l_p$ with $\|\bar{x}_M\|_p = r$, it is clear that $\bar{x} \in \bar{x}_r^{\Uparrow}$. Then, part (b) follows from (a) immediately.

Proof part (c) of (I) in (iii). Let $\bar{x} \in l_p$ with $\|\bar{x}_M\|_p = r$. For any $u \in \bar{x}_r^{\Downarrow}$, there is $\delta > 0$ such that $\|(x + tv)_M\|_p \leq r$, for all $t \in (0, \delta)$. Then, by Lemma 4.2, we have

$$P'_{r\mathbb{C}_M}(\bar{x})(u) = \lim_{t \downarrow 0} \frac{P_{r\mathbb{C}_M}(\bar{x}+tu) - P_{r\mathbb{C}_M}(\bar{x})}{t}$$

$$= \lim_{\substack{t \downarrow 0 \\ t < \delta}} \frac{P_{r\mathbb{C}_M}(\bar{x}+tu) - P_{r\mathbb{C}_M}(\bar{x})}{t}$$

$$= \lim_{\substack{t \downarrow 0 \\ t < \delta}} \frac{(\bar{x}+tu) - \bar{x}}{t}$$

$$= u.$$

This proves (I) of (iii) in this theorem.

Next, we prove (II) in (iii). The proof of (II) is similar to the proof of (II) in part (iii) in Theorem 3.1 in this paper. For an arbitrary given $\bar{x} \in l_p$ with $\|\bar{x}_M\|_p = r$, assume, by the way of contradiction, that $P_{r\mathbb{C}_M}$ is Fréchet differentiable at $\bar{x}$. Then, there is a linear and continuous mapping $A(\bar{x}): l_p \to l_p$, such that

$$\lim_{u \to \bar{x}} \frac{P_{r\mathbb{C}_M}(u) - P_{r\mathbb{C}_M}(\bar{x}) - A(\bar{x})(u - \bar{x})}{\|u - \bar{x}\|_p} = \theta.$$

In particular, in the above limit, we take a directional line segment $u = (1 + \delta)\bar{x}$, for $\delta \downarrow 0$. Since $A(\bar{x})$ is assumed to be linear and continuous, by Lemma 4.2, we have

$$\theta = \lim_{\delta \downarrow 0} \frac{P_{r\mathbb{C}_M}((1+\delta)\bar{x}) - P_{r\mathbb{C}_M}(\bar{x}) - A(\bar{x})((1+\delta)\bar{x} - \bar{x})}{\|(1+\delta)\bar{x} - \bar{x}\|_p}$$

$$= \lim_{\delta \downarrow 0} \frac{\frac{r}{\|((1+\delta)\bar{x})_M\|_p}((1+\delta)\bar{x})_M + ((1+\delta)\bar{x})_{\overline{M}} - \bar{x} - \delta A(\bar{x})(\bar{x})}{\delta\|\bar{x}\|_p}$$

$$= \lim_{\delta \downarrow 0} \frac{\frac{r}{\|((1+\delta)\bar{x})_M\|_p}((1+\delta)\bar{x})_M + ((1+\delta)\bar{x})_{\overline{M}} - \bar{x}_M - \bar{x}_{\overline{M}} - \delta A(\bar{x})(\bar{x})}{\delta\|\bar{x}\|_p}$$

$$= \lim_{\delta \downarrow 0} \frac{\bar{x}_M + ((1+\delta)\bar{x})_{\overline{M}} - \bar{x}_M - \bar{x}_{\overline{M}} - \delta A(\bar{x})(\bar{x})}{\delta\|\bar{x}\|_p}$$

$$= \frac{\bar{x}_{\overline{M}} - A(\bar{x})(\bar{x})}{r}.$$

This implies

$$A(\bar{x})(\bar{x}) = \bar{x}_{\overline{M}}. \tag{4.23}$$

Next, we take an opposite directional line segment $v = (1 - \delta)\bar{x}$, for $\delta \downarrow 0$ with $0 < \delta < 1$. Since $\|\bar{x}_M\|_p = r$, it follows that $(1 - \delta)\bar{x} \in r\mathbb{C}_M$, for any $\delta$ with $0 < \delta < 1$. By Lemma 4.2 and by the assumed linearity of $A(\bar{x})$, we have.

$$\theta = \lim_{\delta \downarrow 0} \frac{P_{r\mathbb{C}_M}((1-\delta)\bar{x}) - P_{r\mathbb{C}_M}(\bar{x}) - A(\bar{x})((1-\delta)\bar{x} - \bar{x})}{\|(1+\delta)\bar{x} - \bar{x}\|_p}$$

$$= \lim_{\delta \downarrow 0} \frac{(1-\delta)\bar{x} - \bar{x} + \delta A(\bar{x})(\bar{x})}{\delta \|\bar{x}\|_p}$$

$$= \lim_{\delta \downarrow 0} \frac{-\delta \bar{x} + \delta A(\bar{x})(\bar{x})}{\delta \|\bar{x}\|_p}$$

$$= \frac{-\bar{x} + A(\bar{x})(\bar{x})}{r}.$$

This implies

$$A(\bar{x})(\bar{x}) = \bar{x}.$$

By $\bar{x}_M \neq \theta$ (by $\|\bar{x}_M\|_p = r$), this contradicts to (4.23), which proves that $P_{r\mathbb{C}_M}$ is not Fréchet differentiable at any point $\bar{x} \in l_p$ with $\|\bar{x}_M\|_p = r$. (II) in (iii) of this theorem is proved. □

In particular, when this given subset $M = \mathbb{N}$, that is $\bar{M} = \emptyset$, then, we have

$$r\mathbb{C}_\mathbb{N} = r\mathbb{B}, \quad \partial(r\mathbb{C}_\mathbb{N}) = r\mathbb{S},$$

and

$$x_\mathbb{N} = x, \quad x_{\bar{\mathbb{N}}} = \theta, \text{ for any } x \in l_p.$$

By Theorem 4.3, we obtain the following results immediately. It is a special case of Theorem 3.1, in which the considered uniformly convex and uniformly smooth Banach space is $l_p$.

**Corollary 4.4.** *For any $r > 0$, the metric projection $P_{r\mathbb{B}}: l_p \to r\mathbb{B}$ has the following differentiability properties.*

(i) $P_{r\mathbb{B}}$ *is strictly Fréchet differentiable on $\mathbb{B}^\circ$ satisfying*

$$\nabla P_{r\mathbb{B}}(\bar{x}) = I_{l_p}, \text{ for every } \bar{x} \in \mathbb{B}^\circ.$$

*That is,*

$$\bar{x} \in \mathbb{B}^\circ \implies \nabla P_{r\mathbb{B}}(\bar{x})(u) = u, \text{ for every } u \in l_p.$$

(ii) $P_{r\mathbb{B}}$ *is Fréchet differentiable at every point $\bar{x} \in l_p \setminus r\mathbb{B}$ such that*

$$\nabla P_{r\mathbb{B}}(\bar{x})(x) = \frac{r}{\|\bar{x}\|_p}\left(u - \frac{\langle J(\bar{x}), u\rangle}{\|\bar{x}\|_p^2}\bar{x}\right), \text{ for every } u \in l_p.$$

*In particular,*

$$\nabla P_{r\mathbb{B}}(\bar{x})(\bar{x}) = \theta, \text{ for every } \bar{x} \in l_p \setminus r\mathbb{B}. \tag{4.7}$$

(iii) *On $r\mathbb{S}$, we have*

    (I) *$P_{r\mathbb{B}}$ is Gâteaux differentiable at every point $\bar{x} \in r\mathbb{S}$ such that*

        (a) $\quad P'_{r\mathbb{B}}(\bar{x})(u) = u - \dfrac{\langle J(\bar{x}), u\rangle}{r^2}\bar{x},\ \ \text{if } u \in \bar{x}_r^{\Uparrow};$

        (b) $\quad P'_{r\mathbb{B}}(\bar{x})(\bar{x}) = \theta;$

        (c) $\quad P'_{r\mathbb{B}}(\bar{x})(u) = u,\ \ \text{if } u \in \bar{x}_r^{\Downarrow}.$

    (II) *$P_{r\mathbb{B}}$ is not Fréchet differentiable at any point in $r\mathbb{S}$, that is,*

$$\nabla P_{r\mathbb{B}}(\bar{x}) \text{ does not exist, for any } \bar{x} \in r\mathbb{S}.$$

## 5. The Fréchet differentiability of the metric projection onto the positive cone in real Banach space $L_p(S)$

Let $(S, \mathcal{A}, \mu)$ be a positive and complete measure space. For any given positive numbers $p$ and $q$ satisfying $1 < p, q < \infty$ and $\dfrac{1}{p} + \dfrac{1}{q} = 1$, in this section, we consider the real uniformly convex and uniformly smooth Banach space $(L_p(S), \|\cdot\|_p)$ with dual space $(L_q(S), \|\cdot\|_q)$. We define the positive cone $K$ in $L_p(S)$ and study the Fréchet differentiability of the metric projection onto $K$.

The normalized duality mapping $J: L_p(S) \to L_q(S)$ has the following representations, for any given $f \in L_p(S)$ with $f \neq \theta$,

$$(Jf)(s) = \frac{|f(s)|^{p-1}\text{sign}(f(s))}{\|f\|_p^{p-2}} = \frac{|f(s)|^{p-2}f(s)}{\|f\|_p^{p-2}}, \text{ for all } s \in S. \tag{5.1}$$

We define a subset $K$ of $L_p(S)$ as follows:

$$K = \{f \in L_p(S): f(s) \geq 0, \text{ for } \mu\text{-almost all } s \in S\}.$$

$K$ is a pointed closed and convex cone in $L_p(S)$ that is called the positive cone of $L_p(S)$.

**Lemma 5.1.** *$K$ has empty interior.*

*Proof.* Let $f \in K$ be arbitrarily given. We show that $f$ is not an interior point of $K$. To this end, for any $\varepsilon > 0$, we can find $\Delta \in \mathcal{A}$ such that

$$\left(\int_\Delta |f(s)|^p \mu(ds)\right)^{\frac{1}{p}} < \frac{\varepsilon}{2} \ \ \text{and}\ \ 0 < \mu(\Delta) < \left(\frac{\varepsilon}{2}\right)^p.$$

Define $g \in L_p(S)$ by

$$g(s) = \begin{cases} f(s), & \text{if } s \notin \Delta, \\ -1, & \text{if } s \in \Delta, \end{cases} \quad \text{for all } s \in S.$$

Then, we have

$$\|f - g\|_p$$

$$= \left( \int_S |f(s) - g(s)|^p \mu(ds) \right)^{\frac{1}{p}}$$

$$= \left( \int_\Delta |f(s) - g(s)|^p \mu(ds) \right)^{\frac{1}{p}}$$

$$\leq \left( \int_\Delta |f(s)|^p \mu(ds) \right)^{\frac{1}{p}} + \left( \int_\Delta |g(s)|^p \mu(ds) \right)^{\frac{1}{p}}$$

$$< \frac{\varepsilon}{2} + (\mu(\Delta))^{\frac{1}{p}}$$

$$< \frac{\varepsilon}{2} + \frac{\varepsilon}{2} = \varepsilon.$$

Since $g \notin K$, this proves that $f$ is not an interior point of $K$. □

In this section, we study the metric projection operator $P_K \colon L_p(S) \to K$. At first, we find the representations of $P_K$.

**Lemma 5.2**. *The metric projection $P_K \colon L_p(S) \to K$ satisfies the following formula.*

$$(P_K f)(s) = \begin{cases} f(s), & \text{if } f(s) \geq 0, \\ 0, & \text{if } f(s) < 0, \end{cases} \quad \text{for any } f \in L_p(S). \tag{5.2}$$

*Proof.* For any given $f \in L_p(S)$, it is clear that (5.2) holds for $f \in K$. So, we suppose that $f \notin K$. By the definition (5.2), if $f \notin K$, then $f - P_K f \neq \theta$. We check the correctness of (5.2) for $f \notin K$. For any $g \in K$, by (5.1) and (5.2), we calculate

$$\langle J(f - P_K f), P_K f - g \rangle$$

$$= \int_S \frac{|(f - P_K f)(s)|^{p-1} \operatorname{sign}((f - P_K f)(s))}{\|f - P_K f\|_p^{p-2}} (P_K f - g)(s) \mu(ds)$$

$$= \int_{f(s) \geq 0} \frac{|(f - P_K f)(s)|^{p-1} \operatorname{sign}((f - P_K f)(s))}{\|f - P_K f\|_p^{p-2}} (P_K f - g)(s) \mu(ds)$$

$$+ \int_{f(s) < 0} \frac{|(f - P_K f)(s)|^{p-1} \operatorname{sign}((f - P_K f)(s))}{\|f - P_K f\|_p^{p-2}} (P_K f - g)(s) \mu(ds)$$

$$= \int_{f(s) \geq 0} \frac{|0|^{p-1} \operatorname{sign}(0)}{\|f - P_K f\|_p^{p-2}} (f - g)(s) \mu(ds) + \int_{f(s) < 0} \frac{|f(s)|^{p-1} \operatorname{sign}(f(s))}{\|f - P_K f\|_p^{p-2}} (0 - g)(s) \mu(ds)$$

$$= \int_{f(s)<0} \frac{|f(s)|^{p-1}(-1)}{\|f-P_Kf\|_p^{p-2}} (-g(s))\mu(ds)$$

$$= \int_{f(s)<0} \frac{|f(s)|^{p-1}}{\|f-P_Kf\|_p^{p-2}} g(s)\mu(ds)$$

$\geq 0$, for all $g \in K$,.

By the basic variational principle of $P_K$, this proves that the correctness of (5.2). □

**Lemma 5.3**. *The metric projection $P_K: L_p(S) \to K$ has the following properties*

(a) $P_K(f) = f$, *for any* $f \in K$;
(b) $P_K(f) = \theta$, *for any* $f \in -K$;
(c) $P_K(f + g) = f + g$, *for any* $f, g \in K$;
(d) *For any* $f \in L_p(S)$,

$$P_K(\lambda f) = \lambda P_K(f), \text{ for any } \lambda \geq 0. \tag{5.3}$$

*Proof.* The proof can be induced by (5.2) immediately. □

By Lemma 5.3, we obtain the following results immediately.

**Proposition 5.4**. *For any $f \in L_p(S)\setminus\{\theta\}$, we have*

$$P'_K(f)(f) = P_K(f). \tag{5.4}$$

*In particular, we have*

$$P'_K(f)(g) = g, \text{ for any } f, g \in K\setminus\{\theta\},$$

$$P'_K(f)(g) = \theta, \text{ for any } f, g \in -K\setminus\{\theta\}.$$

*Proof.* For any $f \in L_p(S)\setminus\{\theta\}$, by (5.3) in Lemma 5.3, we have

$$P'_K(f)(f) = \lim_{t\downarrow 0} \frac{P_K((1+t)f) - P_K(f)}{t}$$

$$= \lim_{t\downarrow 0} \frac{(1+t)P_K(f) - P_K(f)}{t}$$

$$= P_K(f). \qquad \square$$

**Theorem 5.4**. *$P_K$ is not Fréchet differentiable at every point in $L_p(S)$, that is,*

$$\nabla P_K(f) \text{ does not exist, for any } f \in L_p(S).$$

*Proof.* The proof of this theorem is divided into the following two parts:

(i)  $P_K$ is not Fréchet differentiable at every point in $K$, that is,

$$\nabla P_K(f) \text{ does not exist, for any } f \in K. \tag{5.5}$$

(ii) $P_K$ is not Fréchet differentiable at every point in $L_p(S)\setminus K$, that is,

$$\nabla P_K(f) \text{ does not exist, for any } f \in L_p(S)\setminus K. \tag{5.6}$$

Proof of (i). For an arbitrarily given $f \in K$, assume, by the way of contradiction, that $P_K$ is Fréchet differentiable at $f$. Then, there is a linear and continuous mapping $A(f): L_p(S) \to L_p(S)$, such that

$$\lim_{g \to f} \frac{P_K(g) - P_K(f) - A(f)(g-f)}{\|g-f\|_p} = \theta. \tag{5.7}$$

Under the assumption (5.7), we firstly show

$$A(f)(g) = g, \text{ for any } g \in K\setminus\{\theta\}. \tag{5.8}$$

For any $g \in K\setminus\{\theta\}$, in the limit (5.7), we take a directional line segment $f + tg$, for $t \downarrow 0$. Since $A(f)$ is assumed to be linear and continuous, by Lemma 5.3, we have

$$\begin{aligned}
\theta &= \lim_{t \downarrow 0} \frac{P_K(f+tg) - P_K(f) - A(f)(f+tg-f)}{\|f+tg-f\|_p} \\
&= \lim_{t \downarrow 0} \frac{f+tg-f-tA(f)(g)}{t\|g\|_p} \\
&= \lim_{t \downarrow 0} \frac{g - A(f)(g)}{\|g\|_p} \\
&= \frac{g - A(f)(g)}{\|g\|_p}.
\end{aligned}$$

This implies (5.8). Since $f \in K$, there is a positive number $\lambda > 0$ such that

$$\mu\{s \in S: 0 \leq f(s) \leq \lambda\} > 0. \tag{5.9}$$

We take a sequence $\{\Delta_n\}$ of subsets of $\{s \in S: 0 \leq f(s) \leq \lambda\}$ satisfying

$$\Delta_n \in \mathcal{A} \text{ with } \mu(\Delta_n) > 0, \text{ for } n = 1, 2, \ldots$$

and

$$\mu(\Delta_n) \downarrow 0, \text{ as } n \to \infty.$$

For every $n$, we define

$$g_n(s) = \begin{cases} 0, & \text{if } s \notin \Delta_n, \\ 2\lambda, & \text{if } s \in \Delta_n, \end{cases} \text{ for all } s \in S.$$

$g_n \in K\setminus\{\theta\}$, for $n = 1, 2, \ldots$. Then, for every $n$, with respect to the given $f \in K$, we define

$$f_n(s) = \begin{cases} 0, & \text{if } s \notin \Delta_n, \\ f(s), & \text{if } s \in \Delta_n, \end{cases} \text{ for all } s \in S.$$

In the limit (5.7), we take a sequential limit $f - g_n$, for $n \to \infty$. Notice that $f - g_n$ satisfies

$$(f - g_n) = \begin{cases} f(s), & \text{if } s \notin \Delta_n, \\ \leq -\lambda < 0, & \text{if } s \in \Delta_n, \end{cases} \text{ for all } s \in S. \tag{5.10}$$

Since $A(f)$ is assumed to be linear and continuous, by Lemma 5.3, by (5.9) and (5.10), we have

$$\theta = \lim_{n \to \infty} \frac{P_K(f-g_n) - P_K(f) - A(f)(f-g_n-f)}{\|f-g_n-f\|_p}$$

$$= \lim_{n \to \infty} \frac{P_K(f-g_n) - P_K(f) - A(f)(-g_n)}{\|g_n\|_p}$$

$$= \lim_{n \to \infty} \frac{f - f_n - f - A(f)(-g_n)}{\|g_n\|_p}$$

$$= \lim_{n \to \infty} \frac{-f_n + A(f)(g_n)}{\|g_n\|_p}$$

$$= \lim_{n \to \infty} \frac{-f_n + g_n}{\|g_n\|_p}. \tag{5.11}$$

It is clear that $\|g_n\|_p = 2\lambda(\mu(\Delta_n))^{\frac{1}{p}}$. Then, we calculate

$$\left\| \frac{-f_n + g_n}{\|g_n\|_p} \right\|_p^p$$

$$= \frac{1}{(2\lambda)^p \mu(\Delta_n)} \int_S |f(s) - 2\lambda|^p \, \mu(ds)$$

$$= \frac{1}{(2\lambda)^p \mu(\Delta_n)} \int_{\Delta_n} |f(s) - 2\lambda|^p \, \mu(ds)$$

$$\geq \frac{1}{(2\lambda)^p \mu(\Delta_n)} \int_{\Delta_n} \lambda^p \, \mu(ds)$$

$$= \frac{1}{2^p}.$$

This contradicts to (5.11), which proves (5.5).

Proof of (ii). For an arbitrarily given $f \in L_p(S) \backslash K$, assume, by the way of contradiction, that $P_K$ is Fréchet differentiable at $f$. Then, there is a linear and continuous mapping $B(f): L_p(S) \to L_p(S)$, such that

$$\lim_{g \to f} \frac{P_K(g) - P_K(f) - B(f)(g-f)}{\|g-f\|_p} = \theta. \tag{5.12}$$

Since $f \in L_p(S) \backslash K$, there is a positive number $\beta > 0$ such that

$$\mu\{s \in S : 0 > f(s) > -\beta\} > 0.$$

Let $\Delta = \{s \in S: 0 > f(s) > -\beta\}$. Then, we prove the statement: for any $h \in L_p(S)$, if $h$ satisfies the following conditions,

$$h(s) = \begin{cases} 0, & \text{if } s \notin \Delta, \\ \leq 0, & \text{if } s \in \Delta, \end{cases} \text{ for all } s \in S, \tag{5.13}$$

then

$$B(f)(h) = \theta. \tag{5.14}$$

Proof of (5.14). For the given $f \in L_p(S)\setminus K$, for $h \in L_p(S)$ satisfying (5.13), in the limit (5.12), we take a directional line segment $f + th$, for $t \downarrow 0$. Since $B(f)$ is assumed to be linear and continuous, by Lemma 5.3, we have

$$\theta = \lim_{t \downarrow 0} \frac{P_K(f+th) - P_K(f) - B(f)(f+th-f)}{\|f+th-f\|_p}$$

$$= \lim_{t \downarrow 0} \frac{\theta - tB(f)(h)}{t\|h\|_p}$$

$$= \lim_{t \downarrow 0} \frac{-B(f)(h)}{\|h\|_p}$$

$$= \frac{-B(f)(h)}{\|h\|_p}.$$

This proves (5.14). We take a sequence $\{\Delta_m\}$ of subsets of $\Delta$ satisfying

$$\Delta_m \in \mathcal{A} \text{ with } \mu(\Delta_m) > 0, \text{ for } m = 1, 2, \ldots$$

and

$$\mu(\Delta_m) \downarrow 0, \text{ as } m \to \infty.$$

For every $m$, we define

$$h_m(s) = \begin{cases} 0, & \text{if } s \notin \Delta_m, \\ -2f(s), & \text{if } s \in \Delta_m, \end{cases} \text{ for all } s \in S.$$

Since $\Delta_m \subseteq \Delta$, then $-h_m \in L_p(S)$ and $-h_m$ satisfies (5.13), for all $m$. In the limit (5.12), we take a sequential limit $f + h_m$, for $m \to \infty$. Since $B(f)$ is assumed to be linear and continuous, by Lemma 5.3, by the definition of $\Delta$ and by (5.14), we have

$$\theta = \lim_{m \to \infty} \frac{P_K(f+h_m) - P_K(f) - B(f)(f+h_m-f)}{\|f-h_m-f\|_p}$$

$$= \lim_{m \to \infty} \frac{\frac{1}{2}h_m - B(f)(h_m)}{\|h_m\|_p}$$

$$= \lim_{m \to \infty} \frac{\frac{1}{2}h_m + B(f)(-h_m)}{\|h_m\|_p}$$

$$= \frac{1}{2} \lim_{m\to\infty} \frac{h_m}{\|h_m\|_p}.$$

This implies
$$\lim_{m\to\infty} \frac{h_m}{\|h_m\|_p} = \theta.$$

This is a contradiction, which proves (5.6). □